\documentclass[a4paper,11pt]{amsart}

\usepackage{a4wide}
\usepackage[utf8x]{inputenc}
\usepackage[english]{babel}

\usepackage{amsfonts,vmargin,graphicx,amssymb}
\usepackage{psfrag}
\usepackage{amsmath}
\usepackage{xcolor,rotating,epic,eepic}
\usepackage{amsthm}
\usepackage{hyperref}
\usepackage{mathrsfs}
\usepackage{enumerate}
\usepackage{appendix}
\usepackage{url}

\newtheoremstyle{theoreme}
  {10pt}
  {10pt}
  {\itshape}
  {}
  {\bfseries}
  {.}
  {\newline}
  {\thmname{#1}\thmnumber{ #2}\thmnote{ #3}}
  
\newtheoremstyle{lemme}
  {\topsep}
  {\topsep}
  {\itshape}
  {}
  {\bfseries}
  {.}
  {\newline}
  {\thmname{#1}\thmnumber{ #2}\thmnote{ #3}}
  
\newtheoremstyle{remarque}
  {\topsep}
  {\topsep}
  {\normalfont}
  {}
  {\itshape}
  {.}
  {\newline}
  {\thmname{#1}\thmnumber{ #2}\thmnote{ #3}}

\theoremstyle{theoreme}
\newtheorem{theo}{Theorem}[section]

\theoremstyle{lemme}
\newtheorem{lem}[theo]{Lemma}
\newtheorem{prop}[theo]{Proposition}
\newtheorem{dft}[theo]{Definition}

\theoremstyle{remarque}
\newtheorem{rem}[theo]{Remark}

\newcommand{\difx}[1]{\ensuremath{\frac{\partial #1}{\partial x}}}
\newcommand{\difxx}[1]{\ensuremath{\frac{\partial^{2} #1}{\partial x^{2}}}}
\newcommand{\cro}[2]{\ensuremath{\left\langle#1\, ,\, #2\right\rangle}}

\renewcommand{\d}{\ensuremath{\,{\rm d}}}
\newcommand{\sou}[1]{\ensuremath{\underline{#1}}}
\newcommand{\Sun}{\ensuremath{\mathbf{S}^{1}}}
\newcommand{\tend}{\ensuremath{\rightarrow}}
\newcommand{\R}[1]{\ensuremath{\mathbf {R}^{#1}}}
\newcommand{\bZ}{\ensuremath{\mathbf{Z}}}
\newcommand{\om}{\ensuremath{\omega}}
\newcommand{\Som}[3]{\sum_{#1=#2}^{#3}}
\newcommand{\bN}{\ensuremath{\mathbf{N}}}
\newcommand{\mM}{\ensuremath{\mathcal{M}}}
\newcommand{\mS}{\ensuremath{\mathcal{S}}}
\newcommand{\mB}{\ensuremath{\mathcal{B}}}
\newcommand{\mC}{\ensuremath{\mathcal{C}}}
\newcommand{\bD}{\ensuremath{\mathbf{D}}}
\newcommand{\bP}{\ensuremath{\mathbf{P}}}
\newcommand{\bbP}{\ensuremath{\mathbb{P}}}
\newcommand{\mL}{\ensuremath{\mathcal{L}}}
\newcommand{\mN}{\ensuremath{\mathcal{N}}}
\newcommand{\mD}{\ensuremath{\mathcal{D}}}
\newcommand{\mF}{\ensuremath{\mathcal{F}}}
\newcommand{\vphi}{\ensuremath{\varphi}}
\newcommand{\bE}{\ensuremath{\mathbf{E}}}
\newcommand{\bbE}{\ensuremath{\mathbb{E}}}
\newcommand{\mH}{\ensuremath{\mathcal{H}}}
\newcommand{\mR}{\ensuremath{\mathcal{R}}}
\newcommand{\mZ}{\ensuremath{\mathcal{Z}}}
\newcommand{\Om}{\ensuremath{\Omega}}
\newcommand{\eps}{\ensuremath{\varepsilon}}
\newcommand{\N}[1]{\left\|\,#1\,\right\|}
\newcommand{\txi}{\ensuremath{\tilde{\xi}}}
\newcommand{\sig}{\ensuremath{\sigma}}
\newcommand{\ens}[2]{\ensuremath{\left\{#1\,;\,#2\right\}}}
\newcommand{\mU}{\ensuremath{\mathcal{U}}}
\newcommand{\mA}{\ensuremath{\mathcal{A}}}
\newcommand{\Cdc}{\ensuremath{\mathcal{C}^{2}_{c}(\R{})}}

\DeclareMathOperator{\Cov}{Cov}

\title[Quenched limits and fluctuations for rotators in random media]{Quenched limits and fluctuations of the
empirical measure for plane rotators in random media}
\author{Eric Lu{\c c}on}
\date{}
\address{Laboratoire de Probabilit{\'e}s et Mod\`eles Al\'eatoires (CNRS U.M.R. 7599) and  Universit{\'e}
Paris 6
-- Pierre et Marie Curie, U.F.R. Mathematiques, Case 188, 4 place
Jussieu, 75252 Paris cedex 05, France }
\email{eric.lucon\@@etu.upmc.fr}

\begin{document}

\begin{abstract}
The Kuramoto model has been introduced to describe synchronization phenomena observed in groups of cells,
individuals, circuits, etc.
The model consists of $N$ interacting oscillators on the one dimensional sphere $\Sun$, driven by independent
Brownian Motions with constant drift chosen at random. This quenched disorder is chosen independently for each
oscillator according to the same law $\mu$. The behaviour of the system for large $N$ can be understood via
its empirical measure: we prove here the convergence as $N\tend\infty$ of the quenched empirical measure to
the unique solution of coupled McKean-Vlasov equations, under weak assumptions on the disorder $\mu$ and
general hypotheses on the interaction. The main purpose of this work is to address the issue of quenched
fluctuations around this limit, motivated by the dynamical properties of the disordered system for large but
fixed $N$: hence, the main result of this paper is a quenched Central Limit Theorem for the empirical measure.
Whereas we observe a self-averaging for the law of large numbers, this no longer holds for the
corresponding central limit theorem: the trajectories of the fluctuations process are sample-dependent.
\end{abstract}

\maketitle
\section{Introduction}
In this work, we study the fluctuations in the Kuramoto model, which is a particular case of interacting
diffusions with a mean
field Hamiltonian that depends on a random disorder. The Kuramoto model was first introduced in the 70's by
Yoshiki Kuramoto (\cite{cf:Kuramoto}, see also \cite{Acebr'on2005} and references therein) to describe the
phenomenon of synchronization in biological or physical systems. More precisely, the Kuramoto model is a
particular case of a system of $N$ oscillators (considered as elements of the one-dimensional sphere $\Sun:=
\R{}/2\pi\bZ$) solutions to the following SDE:\begin{equation}
\label{eq:introKM}
\d x^{i, N}_{t} = \frac{1}{N} \Som{j}{1}{N}{b(x^{i, N}, x^{j, N}, \om_{j})} \d t + c(x^{i, N}, \om_{i})\d t +
\d B^{i}_{t}, \quad t\in[0, T],\ i=1\dots N,
\end{equation}where $T>0$ is a fixed (but arbitrary) time, $b$ and $c$ are smooth periodic functions. The
Kuramoto case corresponds to a sine interaction ($b(x, y, \om)=K\sin(y-x)$ and $c(x, \om)=\om$). This case
which has the particularity of being rotationally invariant (namely, if $(x^{i, N})_{i}$ is a solution of the
Kuramoto model, $(x^{i, N}+c)_{i}$, $c$ a constant, is also a solution), will be referred to in this work as
the \emph{sine-model}.

The parameter $K>0$ is the coupling strength and $(\om_{j})$ is a sequence of randomly chosen reals (being
i.i.d. realizations of a law $\mu$). The sequence $(\om_{j})_{j}$ is called \emph{disorder} and represents the
fact that the behaviour of each rotator $x^{j, N}$ depends on its own local frequency $\om_{j}$.

Due to the mean field character of \eqref{eq:introKM}, the behaviour of the system can be understood via its
empirical measure $\nu^N$, process with values in $\mM_{1}(\Sun\times\R{})$, that is the set of probability
measures on oscillators and disorder:
\[\forall (\om)\in\R{\bN},\ \forall t\in [0,T],\quad\nu^{N, (\om)}_{t}:= \frac{1}{N}
\Som{j}{1}{N}{\delta_{(x^{j, N}_{t}, \om_{j})}},\]
where $(\om)=(\om_{j})_{j\geq 1}$ is a fixed sequence of disorder in $\R{\bN}$.

\medskip
The purpose of this paper is to address the issue of both convergence and fluctuations of the empirical
measure, as $N\tend\infty$ ; thus the main theorem of this paper (Theorem \ref{theo:fluctquenchedinit})
concerns a Central Limit Theorem in a quenched set-up (namely the quenched fluctuations of $\nu^{N}$ around
its limit).

\medskip
Some heuristic results have been obtained in the physical literature (\cite{Acebr'on2005} and references
therein) concerning the convergence of the empirical measure, as $N\tend\infty$, to a time-dependent measure
$(P_{t}(\d x, \d\om))_{t\in[0,T]}$, whose density w.r.t. Lebesgue measure at time $t$, $q_{t}(x, \om)$ is the
solution of a deterministic non-linear McKean-Vlasov equation (see Eq.\eqref{eq:MKV1L}). It is well understood
(\cite{Acebr'on2005}, \cite{cf:dH}) that crucial features of this equation are captured in the sine-model by
order parameters $r_{t}$ and $\psi_{t}$ defined by: \[r_{t}e^{i\psi_{t}} = \int_{\Sun\times\R{}}e^{ix}q_{t}(x,
\om) \d x\mu(\d\om).\] The quantity $r_{t}$ captures in fact the degree of synchronization of a solution (the
profile $q_t\equiv\frac{1}{2\pi}$ for example corresponds to $r=0$ and represents a total lack of
synchronization) and $\psi_{t}$ identifies the center of synchronization: this is true and rather intuitive
for unimodal profiles. Moreover (\cite{Sakaguchi1988}, see also \cite{cf:dH}, p.118) it turns out that if
$\mu$ is symmetric, all the stationary solutions can be parameterized (up to rotation) by the stationary
version $r$ of $r_{t}$ which must satisfy a fixed point relation $r=\Psi_{K, \mu}(r)$, with $\Psi_{K,
\mu}(\cdot)$ an explicit function such that $\Psi_{K, \mu}(0)=0$. For $K$ small, $r=0$ is the only solution of
such an equation and the system is not synchronized, but for $K$ large, non-trivial solutions appear
(synchronization). In the easiest instances, such a non-trivial solution is unique (in the sense that
$r=\Psi_{K, \mu}(r)$ admits a unique non-zero solution but of course one obtains an infinite number of
solutions by rotation invariance that is $\psi$ can be chosen arbitrarily).

In \cite{cf:dPdH}, Dai Pra and den Hollander have rigorously shown the convergence of the averaged empirical
measure $L^N\in~\mM_{1}(\mC([0,T],\Sun)\times~\R{})$ (probability measure on the whole trajectories and the
disorder): \[L^{N} = \frac{1}{N} \Som{j}{1}{N}{\delta_{(x^{j, N}, \om_{j})}}.\] This convergence of the law of
$L^N$ \emph{under the joint law of the oscillators and the disorder} is shown via an averaged large deviations
principle in the case where $b(x, y, \om)=~K\cdot f(y-~x)$ and $c(x, \om)= g(x,\om)$ for $f$ and $g$ smooth
and \emph{bounded} functions. As a corollary, it is deduced in \cite{cf:dPdH} the convergence of $L_N$ and of
$\nu^N$, via a contraction principle,  in the averaged set-up. In the case of unbounded disorder, the same
proof can be generalized (thesis in preparation) under the following assumption:
\begin{equation}
\label{eq:expmu}
\tag{$H_{\mu}^{A}$}
\forall t>0,\ \int_{\R{}}{e^{t\left|\om\right|}\mu(\d\om)}<\infty.
\end{equation}

One aim of this paper is to obtain the limit of $\nu^{N, (\om)}$ in the \emph{quenched} model, namely for a
\emph{fixed} realization of the disorder $(\om)$. This result can be deduced from the large deviations
estimates in \cite{cf:dPdH}, via a Borel-Cantelli argument, but our result is more direct and works under the
much weaker assumption on $\mu$, $\int_{\R{}}|\om|\mu(\d\om)<\infty$.

\medskip
A crucial aspect of the quenched convergence result, which is a law of large numbers, is that it shows the
\emph{self-averaging} character of this limit: every typical disorder configuration leads as $N\tend\infty$ to
the same evolution equation.

\begin{figure}%
\centering
\includegraphics[width=\textwidth]{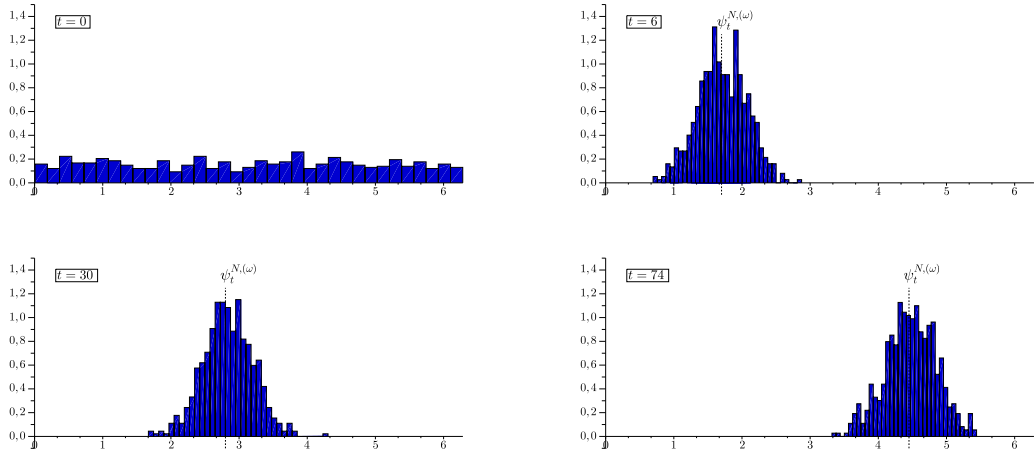}%
\caption{We plot here the evolution of the marginal on $\Sun$ of $\nu^{N, (\om)}$ for $N=600$ oscillators in
the sine-model ($\mu=\frac{1}{2}(\delta_{-1}+ \delta_{1})$, $K=6$). The oscillators are initially chosen
independently and uniformly on $[0, 2\pi]$ independently of the disorder. First the dynamics leads to
synchronization of the oscillators ($t=6$) to a profile which is close to a non-trivial stationary solution of
McKean-Vlasov equation. We then observe that the center $\psi_{t}^{N, (\om)}$ of this density moves to the
right with an approximately constant speed; what is more, this speed of fluctuation turns out to be
sample-dependent (see Fig. \ref{fig:fluctpsiN}).}%
\label{fig:evoldens}%
\end{figure}

\begin{figure}%
\centering
\includegraphics[width=0.8\textwidth]{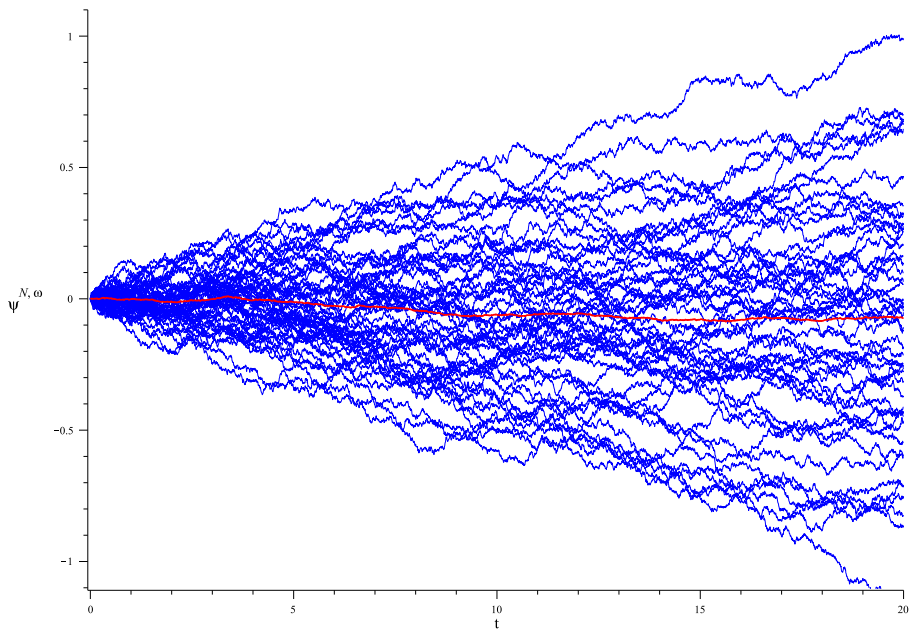}%
\caption{Trajectories of the center of synchronization $\psi^{N, (\om)}$ in the sine-model for different
realizations of the disorder ($\mu=\frac{1}{2}(\delta_{-0.5}+ \delta_{0.5})$, $K=4$, $N=400$). We observe here
the non self-averaging phenomenon: direction and speed of the center depend on the choice of the initial
$N$-sample of the disorder. Moreover, these simulations are compatible with speeds of order $1/\sqrt{N}$.}
\label{fig:fluctpsiN}
\end{figure}

\medskip
However, it seems quite clear even at a superficial level that if we consider the central limit theorem
associated to this convergence, self-averaging does not hold since the fluctuations of the disorder compete
with the dynamical fluctuations. This leads for example to a remarkable phenomenon (pointed out e.g. in
\cite{Balmforth2000} on the basis of numerical simulations): even if the distribution $\mu$ is symmetric, the
fluctuations of a fixed chosen sample of the disorder makes it \emph{not symmetric} and thus the center of the
synchronization of the system \emph{slowly} (i.e. with a speed of order $1/\sqrt{N}$) rotates in one direction
and with a speed that depends on the sample of the disorder (Fig. \ref{fig:evoldens} and \ref{fig:fluctpsiN}).
This non-self averaging phenomenon can be tackled in the sine-model by computing the finite-size order
parameters (Fig. \ref{fig:fluctpsiN}): \begin{equation}\label{eq:rNpsiN}r^{N, (\om)}_{t} e^{i\psi^{N,
(\om)}_{t}} = \frac{1}{N} \Som{j}{1}{N}{e^{i x_{t}^{j, N}}} = \cro{\nu_{t}^{N, (\om)}}{e^{ix}}.\end{equation}

\medskip
As a step toward understanding this non self-averaging phenomenon, the second and main goal of this paper is
to
establish a fluctuations result around the McKean-Vlasov limit in the quenched set-up (see Theorem
\ref{theo:fluctquenchedinit}). A central limit theorem
for the averaged model is addressed in \cite{cf:dPdH}, applying techniques introduced by Bolthausen
\cite{Bolthausen1986}. A fluctuations theorem may also be found in
\cite{Collet2010} for a model of social interaction in an averaged set-up. Here we
prove convergence in law of the \emph{quenched} fluctuations process
\[\eta^{N, (\om)}:= \sqrt{N} \left(\nu^{N, (\om)} - P\right),\] seen as a continuous process in the Schwartz
space $\mS'$ of tempered distributions on $\Sun\times\R{}$ to the solution of an Ornstein-Uhlenbeck process.
The quenched convergence is here understood as a weak convergence \emph{in law w.r.t. the disorder} and is
more technically involved than the convergence in the averaged system. The main techniques we exploit have
been introduced by Fernandez and M\'el\'eard \cite{cf:FernMeleard} and Hitsuda and Mitoma \cite{cf:Hitsuda},
who studied similar fluctuations in the case without disorder. In \cite{BudDupFisch2010}, A Large Deviation
Principle is also proved. We refer to Section \ref{sec:fluct} for
detailed definitions.
\medskip

While numerical computations of the trajectories of the limit process of fluctuations clearly show a non
self-averaging phenomenon, the dynamical properties of the fluctuations process that we find are not
completely understood so far. Progress in this direction requires a good understanding of the spectral
properties of the linearized operator of McKean-Vlasov equation around its non-trivial stationary solution;
the stability of the non-synchronized solution $q\equiv\frac{1}{2\pi}$ has been treated by Strogatz and
Mirollo in \cite{Strogatz1991}. In the particular case without disorder, spectral properties of the evolution
operator linearized around the non-trivial stationary solution are obtained in \cite{BertGiacPak}, but the
case with a general distribution $\mu$ needs further investigations.

This work is organized as follows: Section \ref{sec:modtheo} introduces the model and the main results.
Section \ref{sec:convps} focuses on the quenched convergence of $\nu_N$. In Section \ref{sec:fluct}, the
quenched Central Limit Theorem is proved. The last section \ref{sec:prooffluctrpsi} applies the fluctuations
result to the behaviour of the order parameters in the sine-model.

\section{Notations and main results}
\label{sec:modtheo}
\subsection{Notations}
\label{subsec:not}
\begin{itemize}
	\item if $X$ is a metric space, $\mB_{X}$ will be its Borel $\sigma$-field,
	\item $\mC_{b}(X)$ (resp. $\mC_{b}^{p}(X)$, $p=1, \dots, \infty$), the set of bounded continuous
functions (resp. bounded continuous with bounded continuous derivatives up to order $p$) on $X$, ($X$ will be
often $\Sun\times \R{}$),
	\item $\mC_{c}(X)$ (resp. $\mC_{c}^{p}(X)$, $p=1, \dots, \infty$), the set of continuous functions
with compact support (resp. continuous with compact support with continuous derivatives up to order $p$) on
$X$,
	\item $\bD([0,T], X)$, the set of right-continuous with left limits functions with values on $X$,
endowed with the Skorokhod topology,
	\item $\mM_{1}(Y)$, the set of probability measures on $Y$ ($Y$ topological space, with a regular
$\sigma$-field $\mB$),
	\item $\mM_{F}(Y)$, the set of finite measures on $Y$,
	\item $(\mM_{1}(Y),w)$: $\mM_{1}(Y)$ endowed with the topology of weak convergence, namely the
coarsest topology on $\mM_{1}(Y)$ such that the evaluations $\nu\mapsto \int{f\d\nu}$ are continuous, where
$f$ are bounded continuous,
	\item $(\mM_{1}(Y),v)$: $\mM_{1}(Y)$ endowed with the topology of vague convergence, namely the
coarsest topology on $\mM_{1}(Y)$ such that the evaluations $\nu\mapsto \int{f\d\nu}$ are continuous, where
$f$ are continuous with compact support.
\end{itemize}
We will use $C$ as a constant which may change from a line to another.
\subsection{The model}
We consider the solutions of the following system of SDEs:

for $i=1,\dots, N$, for $T>0$, for all $t\leq T$,
\begin{equation}
\label{eq:edsbc}
x^{i,N}_{t} = \xi^{i} + \frac{1}{N}\Som{j}{1}{N}{\int_{0}^{t}{b(x_{s}^{i, N}, x_{s}^{j, N}, \om_{j})\d s}} +
\int_{0}^{t}{c(x_{s}^{i, N},\om_{i})\d s} + B^{i}_t,
\end{equation}
where the initial conditions $\xi^{i}$ are independent and identically distributed with law $\lambda$, and
independent of the Brownian motion $(B) = (B^i)_{i\geq 1}$, and where $b$ (resp. $c$) is a smooth function,
$2\pi$-periodic w.r.t. the two first (resp. first) variables. The disorder $(\om)= (\om_{i})_{i\geq 1}$ is a
realization of i.i.d. random variables with law $\mu$.
\begin{rem}
The assumption that the random variables $(\om_{i})$ are independent will not always be necessary and will be
weakened when possible.
\end{rem}
Instead of considering $x^{i,N}$ as elements of $\R{}$, we will consider their projection on $\Sun$. For
simplicity, we will keep the same notation $x^{i,N}$ for this projection\footnote{See Remarks
\ref{rem:noncomp1} and \ref{rem:noncomp2} for possible generalizations to the non-compact case.}.

We introduce the empirical measure $\nu^{N}$ (on the trajectories and disorder):
\begin{dft}
For all $t\leq T$, for a fixed trajectory $(x^1,\dots, x^N) \in \mathcal{C}([0, T], (\Sun)^{N})$ and a fixed
sequence of disorder $(\om)$, we define an element of $\mathcal{M}_{1}(\Sun\times \R{})$ by:
\begin{equation*}
\label{eq:muN1}
\nu^{N, (\om)}_{t} = \frac{1}{N} \Som{i}{1}{N}{\delta_{(x^{i,N}_{t}, \om_{i})}}.
\end{equation*}
\end{dft}

Finally, we introduce the fluctuations process $\eta^{N, (\om)}$ of $\nu^{N, (\om)}$ around its limit $P$ (see
Th.~\ref{prop:convps}):
\begin{dft}
For all $t\leq T$, for fixed $(\om)\in \R{\bN}$, we define:
\begin{equation*}
\label{eq:etaN}
\eta^{N, (\om)}_{t} = \sqrt{N}\left(\nu^{N, (\om)}_{t} - P_t\right).
\end{equation*}
\end{dft}

Throughout this article, we will denote as $\bP$ the law of the sequence of Brownian Motions and as
$\mathbb{P}$ the law of the sequence of the disorder. The corresponding expectations will be denoted as $\bE$
and $\mathbb{E}$ respectively.
\subsection{Main results}
\subsubsection{Quenched convergence of the empirical measure}
\label{subsubsec:quenEM}
In \cite{cf:dPdH}, Dai Pra and den Hollander are interested in the \emph{averaged} model, i.e. in the
convergence in law of the empirical measure \emph{under the joint law of both oscillators and disorder}.
The model studied here, which is more interesting as far as the biological applications are concerned is
\emph{quenched}: for a fixed realization of the disorder $(\om)$, do we have the convergence of the empirical
measure? Moreover the convergence is shown under weaker assumptions on the moments of the disorder.

We consider here the general case where $b(x, y, \om)$ is bounded, Lipschitz-continuous, and $2\pi$-periodic
w.r.t. the two first variables. $c$ is assumed to be Lipschitz-continuous w.r.t. its first variable, but not
necessarily bounded (see the sine-model, where $c(x, \om)= \om$). We also suppose that the function
$\om\mapsto S(\om) := \sup_{x\in\Sun} |c(x, \om)|$ is continuous (this is in particular true if $c$ is
uniformly continuous w.r.t. to both variables $(x, \om)$, and obvious in the sine-model where $S(\om)=|\om|$).
The Lipschitz bounds for $b$ and $c$ are supposed to be uniform in $\om$.
 
The disorder $(\om)$, is assumed to be a sequence of identically distributed random variables (but not
necessarily independent), such that the law of each $\om_{i}$ is $\mu$. We suppose that the sequence $(\om)$
satisfies the following property: for $\mathbb{P}$-almost every sequence $(\om)$,
\begin{equation}
\label{eq:Hmu}
\tag{$H_{\mu}^{Q}$}
\frac{1}{N}\Som{i}{1}{N}{\sup_{x\in\Sun}|c(x, \om_{i})|} \tend_{N\tend\infty} \int{\sup_{x\in\Sun}|c(x,
\om)|\mu(\d\om)}<\infty.
\end{equation}

We make the following hypothesis on the initial empirical measure:

\begin{equation}
\label{eq:H3nu}
\tag{$H_{0}$}
\nu^{N, (\om)}_{0} \stackrel{N\tend\infty}{\longrightarrow} \nu_{0}, \quad \textit{in law, in
$(\mM_{1}(\Sun\times\R{}), w)$}.
\end{equation}

\begin{rem}
\begin{enumerate}
	\item The required hypotheses about the disorder and the initial conditions are weaker than for the
large deviation principle:
\begin{itemize}
	\item the (identically distributed) variables $(\om_{i})$ need not be independent: we simply need a
convergence (similar to a law of large numbers) only concerning the function $S$,
	\item Condition \eqref{eq:Hmu} is weaker than \eqref{eq:expmu} on page \pageref{eq:expmu}; for the
sine-model, \eqref{eq:Hmu} reduces to $\int{|\om|\mu(\d\om)}<\infty$,
	\item the initial values need not be independent, we only assume a convergence of the empirical
measure.
\end{itemize}

\item The hypothesis \eqref{eq:Hmu} is verified, for example, in the case of i.i.d. random variables, or in
the case of an ergodic stationary Markov process.
\item Under \eqref{eq:H3nu}, the second marginal of $\nu_{0}$ is $\mu$.
\end{enumerate}
\end{rem}

In Section \ref{sec:convps}, we show the following:

\begin{theo}\label{prop:convps}
Under the hypothesis \eqref{eq:H3nu} and \eqref{eq:Hmu}, for $\bbP$-almost every sequence $(\om_{i})$, the
random variable $\nu^{N, (\om)}$ converges in law to $P$, in the space $\bD([0, T], (\mM_{1}(\Sun\times\R{}),
w))$, where $P$ is the only solution of the following weak equation (for every $f$ continuous bounded on
$\Sun\times\R{}$, twice differentiable, with bounded derivatives):
\begin{equation}
\cro{P_{t}}{f} = \cro{\nu_{0}}{f} + \frac{1}{2} \int_{0}^{t}{\d s \cro{P_{s}}{f''}} + \int_{0}^{t}{\d s
\cro{P_{s}}{f' (b[\cdot , P_{s}] + c)}},
\label{eq:weakMKV}
\end{equation}
where \[b[x, m]= \int{b(x,y,\pi)m(\d y, \d\pi)}.\]
Moreover, with the same hypotheses, the law of $\nu^N$ under the joint law of the oscillators and disorder
(averaged model) converges weakly to $P$ as well.
\end{theo}

\begin{rem}
An easy calculation shows that $P$ can be considered as a weak solution to the family of coupled McKean-Vlasov
equations (see \cite{cf:dPdH}):
\begin{enumerate}
	\item $P$ can be written as $P(\d x, \d\om) = \mu(\om) P^{\om}(\d x)$,
	\item if we define $q_{t}^{\om}$ through $P_{t}(\d x, \d\om) = \mu(\om) q_{t}^{\om}(\d x)$,
$q_{t}^{\om}$ is the unique weak solution of the McKean-Vlasov equation:
\begin{align}
\frac{\d }{\d t} q_{t}^{\om} &= \mL^{\om} q_{t}^{\om},\quad q_{0}^{\om} =\lambda.\label{eq:MKV1L}
\end{align}
where, $\mL^{\om}$ is the following differential operator:
\begin{equation}
\label{eq:Lom}
\mL^{\om}q_{t}^{\om} = - \difx{}\left[\left(\int_{\R{}}{b(x, y, \pi) q_{t}^{\pi}(\d y)\mu(\d\pi)}+ c(x,
\om)\right)q_{t}^{\om}\right] + \frac{1}{2} \difxx{}q_{t}^{\om}.
\end{equation}
\end{enumerate}
We insist on the fact that Eq. \eqref{eq:MKV1L} is indeed a (possibly) infinite system of coupled non-linear
PDEs. To fix ideas, one may consider the simple case where $\mu=\frac{1}{2}(\delta_{-1}+\delta_{1})$. Then
\eqref{eq:MKV1L} reduces to two equations (one for $+1$, the other for $-1$) which are coupled via the
averaged measure $\frac{1}{2}(q_{t}^{+1}+ q_{t}^{-1})$. But for more general situations ($\mu=\mN(0, 1)$ say)
this would consist of an infinite number of coupled equations.
\end{rem}
\mbox{}

\begin{rem}[(Generalization to the non compact case)]
\label{rem:noncomp1}
 The assumption that the state variables are in $\Sun$, although motivated by the Kuramoto model, is not
absolutely essential: Theorem \ref{prop:convps} still holds in the non-compact case (e.g. when $\Sun$ is
replaced by $\R{d}$), under the additional assumptions of boundedness of $x\mapsto |c(x, \om)|$
and the first finite moment of the initial condition: $\int |x| \lambda(\d x)<\infty$.
\end{rem}
\medskip
We know turn to the statement of the main result of the paper: Theorem \ref{theo:fluctquenchedinit}.
\subsubsection{Quenched fluctuations of the empirical measure}
Theorem \ref{prop:convps} says that for $\bbP$-almost every realization $(\om)$ of the disorder, we have the
convergence of $\nu^{N,(\om)}$ towards $P$, which is a law of large numbers. We are now interested in the
corresponding Central Limit Theorem associated to this convergence, namely, for a
\emph{fixed} realization of the disorder $(\om)$, in the asymptotic behaviour, as $N\tend\infty$ of the
fluctuations field $\eta^{N, (\om)}$ taking values in the set of signed measures:
\[\forall t\in[0, T],\ \eta_{t}^{N, (\om)}:= \sqrt{N} \left(\nu_{t}^{N, (\om)} - P_{t}\right).\]

In the case with no disorder, such fluctuations have already been studied by numerous authors (eg. Sznitman
\cite{cf:Sznit84}, Fernandez-M\'el\'eard \cite{cf:FernMeleard}, Hitsuda-Mitoma \cite{cf:Hitsuda}). More
particularly, Fernandez and M\'el\'eard show the convergence of the fluctuations field in an appropriate
Sobolev space to an Ornstein-Uhlenbeck process.

Here, we are interested in the \emph{quenched} fluctuations, in the sense that the fluctuations are studied
for fixed realizations of the disorder. We will prove a weak convergence of the law of the process $\eta^{N,
(\om)}$, \emph{in law w.r.t. the disorder}.

In addition to the hypothesis made in \S \ref{subsubsec:quenEM}, we make the following assumptions about $b$
and $c$ (where $\mD_{p}$ is the set of all differential operators of the form
$\partial_{u^{k}}\partial_{\pi^{l}}$ with $k+l\leq p$):

\begin{equation}
\label{eq:Hbc}
\tag{$H_{b,c}$}\left\{\begin{array}{c}
b \in \mC_{b}^{\infty}(\Sun\times \R{}),\quad c \in \mC^{\infty}(\Sun\times \R{}),\vspace{10pt}\\
\displaystyle\exists \alpha >0, \sup_{D\in \mD_{6}}\int_{\R{}}{\frac{\sup_{u\in\Sun}|D c(u,
\pi)|^{2}}{1+|\pi|^{2\alpha}}\d\pi} <\infty,
\end{array}\right.
\end{equation}
Furthermore, we make the following assumption about the law of the disorder ($\alpha$ is defined in
\eqref{eq:Hbc}):
\begin{equation}
\label{eq:Hmufluct}
\tag{$H_{\mu}^{\mF}$}
\text{the $(\om_{j})$ are i.i.d. and}\quad \int_{\R{}}{|\om|^{4\alpha}\mu(\d\om)}<\infty.
\end{equation}
\begin{rem}
The regularity hypothesis about $b$ and $c$ can be weakened (namely $b \in \mC_{b}^{n}(\Sun\times
\R{})$ and $c \in \mC^{m}(\Sun\times \R{})$ for sufficiently large $n$ and $m$) but we have kept $m=n=\infty$
for the sake of clarity.
\end{rem}
\begin{rem}
In the case of the sine-model, Hypothesis \eqref{eq:Hbc} is satisfied with $\alpha=2$ for example.
\end{rem}

In order to state the fluctuations theorem, we need some further notations: for all $s\leq T$, let $\mL_{s}$
be the second order differential operator defined by 
\[\mL_{s}(\vphi)(y, \pi):=  \frac{1}{2} \vphi''(y, \pi) + \vphi'(y, \pi)(b[y, P_{s}] + c(y, \pi)) +
\cro{P_{s}}{\vphi'(\cdot, \cdot)b(\cdot, y, \pi)}.\]
Let $W$ the Gaussian process with covariance: \begin{equation}
\label{eq:Wmart}
\bE(W_{t}(\vphi_{1})W_{s}(\vphi_{2})) = \int_{0}^{s\wedge t}{\cro{P_{u}}{\vphi'_{1}\vphi'_{2}}\d u}.
\end{equation}
For all $\vphi_{1}, \vphi_{2}$ bounded and continuous on $\Sun\times \R{}$, let \begin{align}
 \Gamma_{1}(\vphi_{1}, \vphi_{2}) &= \int_{\R{}}{\Cov_{\lambda}\left(\vphi_{1}(\cdot\ ,\om), \vphi_{2}(\cdot\
, \om)\right)\mu(\d\om)},\label{eq:covC1}\\ &= \int_{\Sun\times\R{}}{{\left(\vphi_{1} -
\int_{\Sun}{\vphi_{1}(\cdot\ , \om)\d\lambda}\right)\left(\vphi_{2} - \int_{\Sun}{\vphi_{2}(\cdot\ ,
\om)\d\lambda}\right)\lambda(\d x)}\mu(\d\om)},\nonumber
\end{align} and \begin{align}
 \Gamma_{2}(\vphi_{1}, \vphi_{2}) &=\Cov_{\mu}\left(\int_{\Sun}\vphi_{1}\d\lambda,
\int_{\Sun}\vphi_{2}\d\lambda\right),\label{eq:covC2}\\&=\int_{\R{}}{\left(\int_{\Sun}{\vphi_{1}\d\lambda}-
\int_{\Sun\times\R{}}{\vphi_{1}\d\lambda\d\mu}\right)}\left(\int_{\Sun}{\vphi_{2}\d\lambda}-
\int_{\Sun\times\R{}}{\vphi_{2} \d\lambda \d\mu}\right)\d\mu.\nonumber
\end{align}

For fixed $(\om)$, we may consider $\mH_{N}(\om)$, the law of the process $\eta^{N, (\om)}$; $\mH_{N}(\om)$
belongs to $\mM_{1}(\mC([0, T], \mS'))$, where $\mS'$ is the usual Schwartz space of tempered distributions on
$\Sun\times\R{}$. We are here interested in the law of the random variable $(\om) \mapsto \mH_{N}(\om)$ which
is hence an element of $\mM_{1}(\mM_{1}(\mC([0, T], \mS')))$.

The main theorem (which is proved in Section \ref{sec:fluct}) is the following:
\begin{theo}[(Fluctuations in the quenched model)]
\label{theo:fluctquenchedinit}
Under \eqref{eq:Hmufluct}, \eqref{eq:Hbc}, the sequence $(\om)\mapsto \mH_{N}(\om)$
converges in law to the random variable $\om\mapsto \mH(\om)$, where $\mH(\om)$ is the law of the solution to
the Ornstein-Uhlenbeck process $\eta^{\om}$ solution in $\mS'$ of the following equation:
\begin{equation}
\label{eq:OUeta}
\eta_{t}^{\om} = X(\om) + \int_{0}^{t}{\mL^{*}_{s}\eta_{s}^{\om}\d s} + W_{t},
\end{equation}
where, $\mL^{*}_{s}$ is the formal adjoint operator of $\mL_{s}$ and for all fixed $\om$, $X(\om)$ is
a non-centered Gaussian process with covariance $\Gamma_{1}$ and with mean value $C(\om)$. As a random
variable in $\om$, $\om\mapsto C(\om)$ is a Gaussian process with covariance $\Gamma_{2}$. Moreover, $W$ is
independent on the initial value $X$.
\end{theo}
\begin{rem}
In the evolution \eqref{eq:OUeta}, the linear operator $\mL_{s}^{*}$ is deterministic ; the only dependence
in $\om$ lies in the initial condition $X(\om)$, through its non trivial means $C(\om)$. However,
numerical simulations of trajectories of $\eta^{\om}$ (see Fig. \ref{fig:evoleta}) clearly show a non
self-averaging phenomenon, analogous to the one observed in Fig \ref{fig:fluctpsiN}: $\eta^{\om}_{t}$ not only
depends on $\om$ through its initial condition $X(\om)$, but also for all positive time $t>0$.

Understanding how the deterministic operator $\mL_{s}^{*}$ propagates the initial dependence in $\om$ on the
whole trajectory is an intriguing question which requires further investigations (work in progress). In that
sense, one would like to have a precise understanding of the spectral properties of $\mL_{s}^{*}$, which
appears to be deeply linked to the differential operator in McKean-Vlasov
equation \eqref{eq:Lom} linearized around its non-trivial stationary solution.
\end{rem}
\mbox{}

\begin{rem}[(Generalization to the non-compact case)]
\label{rem:noncomp2}
 As in Remark \ref{rem:noncomp1}, it is possible to extend Theorem \ref{theo:fluctquenchedinit} to the
(analogous but more technical) case where $\Sun$ is replaced by $\R{d}$. To this purpose, one has to
introduce an additional weight $(1+|x|^\alpha)^{-1}$ in the definition of the Sobolev norms in Section
\ref{sec:fluct} and to suppose appropriate hypothesis concerning the first moments of the initial condition
$\lambda$ ($\int |x|^{\beta}\lambda(\d x)<\infty$ for a sufficiently large $\beta$).
\end{rem}

\subsubsection{Fluctuations of the order parameters in the Kuramoto model}
\label{sec:fluctrpsi}
For given $N\geq 1$, $t\in[0,T]$ and disorder $(\om)\in\R{\bN}$, let $r_{t}^{N, (\om)}>0$ and $\zeta_{t}^{N,
(\om)}\in\Sun$ such that
\begin{equation*}
r_{t}^{N, (\om)} \zeta_{t}^{N, (\om)} = \frac{1}{N} \Som{j}{1}{N}{e^{i x_{t}^{j, N}}} = \cro{\nu_{t}^{N,
(\om)}}{e^{ix}}.
\label{eq:rpsiN}
\end{equation*}

\begin{prop}[(Convergence and fluctuations for $r_{t}^{N, (\om)}$)] We have the following:
\label{prop:convfluctr}
\begin{enumerate}
	\item Convergence of $r_{t}^{N, (\om)}$:
For $\bbP$-almost every realization of the disorder, $r^{N, (\om)}$ converges in law in $\mC([0, T], \R{})$,
to $r$ defined by 
\begin{equation*}
t\in[0,T]\mapsto r_{t}:= \left(\cro{P_{t}}{\cos(\cdot)}^{2}+\cro{P_{t}}{\sin(\cdot)}^{2}\right)^{\frac{1}{2}}.
\label{eq:rt}
\end{equation*}
\item If $r_{0}>0$ then 
\begin{equation}
\forall t\in[0,T],\quad r_{t}>0.
\label{eq:rtpos}
\tag{$H_{r}$}
\end{equation}
\item Fluctuations of $r_{t}^{N, (\om)}$ around its limit:
Let \begin{equation*}
t\mapsto \mR_{t}^{N,(\om)}:= \sqrt{N}\left(r_{t}^{N, (\om)} - r_{t}\right)
\label{eq:etar}
\end{equation*} be the fluctuations process. For fixed disorder  $(\om)$, let $\mathfrak{R}^{N, (\om)}\in
\mM_{1}(\mC([0, T], \R{})) $ be the law of $\mR^{N, (\om)}$.
Then, under \eqref{eq:rtpos}, the random variable $(\om)\mapsto \mathfrak{R}^{N, (\om)}$ converges in law to
the random variable $\om\mapsto \mathfrak{R}^{\om}$, where $\mathfrak{R}^{\om}$ is the law of
$\mR^{\om}:=\frac{1}{r}\left(\cro{P}{\cos(\cdot)}\cdot\cro{\eta^{\om}}{\cos(\cdot)}+\cro{P}{\sin(\cdot)}
\cdot\cro{\eta^{\om}}{\sin(\cdot)}\right)$.
\end{enumerate}
\end{prop}

\begin{figure}%
\centering
\includegraphics[width=\textwidth]{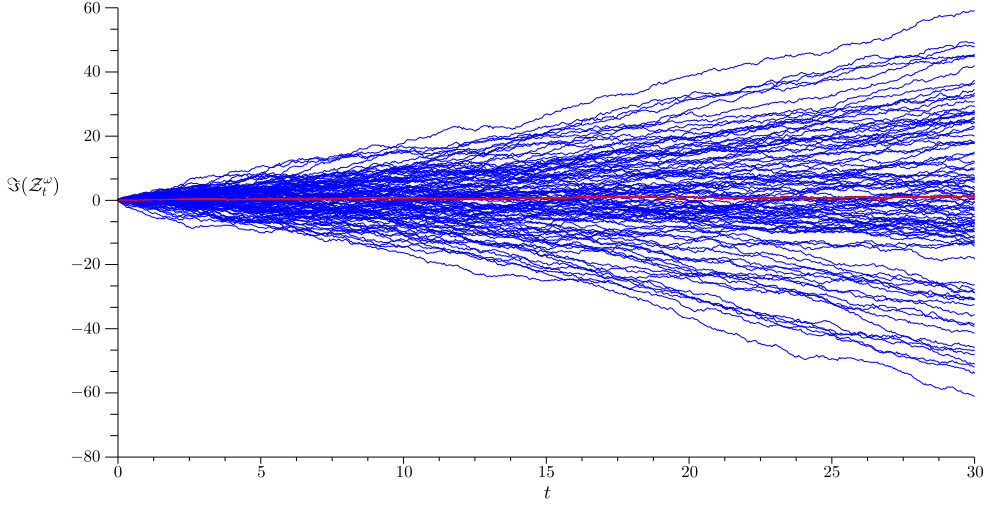}%
\caption{We plot here the evolution of the imaginary part of the process $\mZ^{\om}$, for different
realizations of $\om$; the trajectories are sample-dependent, as in Fig \ref{fig:fluctpsiN}.}%
\label{fig:evoleta}%
\end{figure}

\begin{rem}
\label{rem:Rtp}
In simpler terms, this double convergence in law corresponds for example to the convergence in law of the
corresponding characteristic functions (since the tightness is a direct consequence of the tightness of the
process $\eta$); i.e. for $t_{1}, \dots, t_{p}\in[0, T]$ ($p\geq 1$) the characteristic function of
$(\mR_{t_{1}}^{N, (\om)},\dots, \mR_{t_{p}}^{N, (\om)})$ for fixed $(\om)$ converges in law, as a random
variable in $(\om)$, to the random characteristic function of $(\mR^{\om}_{t_{1}}, \dots, \mR^{\om}_{t_{p}})$.
\end{rem}
\medskip
\begin{prop}[(Convergence and fluctuations for $\zeta^{N, (\om)}$)]We have the following:
\label{prop:convfluctpsi}
\begin{enumerate}
	\item Convergence of $\zeta^{N, (\om)}$: Under \eqref{eq:rtpos}, for $\bbP$-almost every realization
of the disorder $(\om)$, $\zeta^{N, (\om)}$ converges in law to $\zeta: t\in[0, T]\mapsto \zeta_{t}:=
\frac{\cro{P_{t}}{e^{ix}}}{r_{t}}$,
\item Fluctuations of $\zeta^{N, (\om)}$ around its limit:
Let \begin{equation*}
t\mapsto \mZ_{t}^{N,(\om)}:= \sqrt{N}\left(\zeta_{t}^{N, (\om)} - \zeta_{t}\right)
\label{eq:etapsi}
\end{equation*} be the fluctuations process. For fixed disorder $(\om)$, let $\mathfrak{Z}^{N, (\om)}\in
\mM_{1}(\mC([0, T], \R{})) $ be the law of $\mZ^{N, (\om)}$. Then, under \eqref{eq:rtpos}, the random variable
$(\om)\mapsto \mathfrak{Z}^{N, (\om)}$ converges in law to the random variable $\om\mapsto
\mathfrak{Z}^{\om}$, where $\mathfrak{Z}^{\om}$ is the law of
$\mZ^{\om}:=\frac{1}{r^2}\left(r\cro{\eta^{\om}}{\cos(\cdot)}+\cro{P}{e^{ix}}\mR^{\om}\right)$.
\end{enumerate}
\end{prop}
In the sine-model, we have $\zeta^{N, (\om)}= e^{i\psi^{N, (\om)}}$ where $\psi^{N, (\om)}$ is defined in
\eqref{eq:rNpsiN} and is plotted in Fig. \ref{fig:fluctpsiN}. Some trajectories of the process $\mZ^{\om}$
are plotted in Fig. \ref{fig:evoleta}.

This fluctuations result is proved in Section \ref{sec:prooffluctrpsi}.

\section{Proof of the quenched convergence result}
\label{sec:convps}
In this section we prove Theorem \ref{prop:convps}. Reformulating \eqref{eq:edsbc} in terms of $\nu^{N,
(\om)}$, we have:
\begin{equation}
\label{eq:diffbc}
\forall i=1\dots N,\ \forall t\in[0, T],\ x^{i,N}_{t} = \xi^{i} +\int_{0}^{t}{b[x_{s}^{i, N}, \nu_{s}^{N}]\d
s} + \int_{0}^{t}{c(x_{s}^{i, N},\om_{i})\d s} + B^{i}_t,
\end{equation} where we recall that $b[x, m]:= \int{b(x, y, \pi)m(\d y, \d\pi)}$.

The idea of the proof of Theorem \ref{prop:convps} is the following: we show the tightness of the sequence
$(\nu^{N, (\om)})$ firstly in $\bD([0, T], (\mM_{F}, v))$ (recall Notations in \S \ref{subsec:not}), which is
quite simple since $\mC_{c}(\Sun\times\R{})$ is separable and by an argument of boundedness of the second
marginal of any accumulation point, thanks to \eqref{eq:Hmu}, we show the tightness in $\bD([0, T], (\mM_{F},
w))$.
The proof is complete when we prove the uniqueness of any accumulation point.

\subsection{Proof of the tightness result}
We want to show successively:

\begin{enumerate}
	\item Tightness of $\mL(\nu^{N, (\om)})$ in $\bD([0, T], (\mM_{F}, v))$,
	\item Equation verified by any accumulation point,
	\item Characterization of the marginals of any limit,
	\item Convergence in $\bD([0, T], (\mM_{F}, w))$.
\end{enumerate}

\subsubsection{Equation verified by $\nu^{N, (\om)}$}
For $f\in\mC^{2}_{b}(\Sun\times \R{})$, we denote by $f'$, $f''$ the first and second derivative of $f$ with
respect to the first variable. Moreover, if $m\in\mM_{1}(\Sun\times\R{})$, then $\cro{m}{f}$ stands for
$\int_{\Sun\times\R{}}f(x, \pi) m(\d x, \d\pi)$.

Applying Ito's formula to \eqref{eq:diffbc}, we get, for all $f\in\mC^{2}_{b}(\Sun\times \R{})$,
\begin{align*}
	\cro{\nu^{N, (\om)}_{t}}{f} &= \cro{\nu^{N, (\om)}_{0}}{f} + \frac{1}{2} \int_{0}^{t}{\d s
\cro{\nu^{N, (\om)}_{s}}{f''}}\nonumber\\ &+ \int_{0}^{t}{\d s \cro{\nu^{N, (\om)}_{s}}{f'\cdot(b[\cdot,
\nu^{N, (\om)}_{s}] + c)}}+ M_{N, f}(t),\label{eq:nufdis}
\end{align*}
where $M_{N, f}(t):= \frac{1}{N} \Som{j}{1}{N}{\int_{0}^{t}{f'(x_{j}^{N, (\om)}, \om_{j})\d B_{j}(s)}}$ is a
martingale ($f'$ bounded).

\subsubsection{Tightness of $\mL(\nu^{N, (\om)})$ in $\bD([0, T], (\mM_{F}, v))$}
$\mC_{c}(\Sun\times\R{})$ is separable: let $(f_{k})_{k\geq1}$ (elements of $\mC^{\infty}(\Sun\times\R{})$) a
dense sequence in $\mC_{c}(\Sun\times\R{})$, and let $f_{0}\equiv 1$. We define $\Om:= \bD([0, T], (\mM_{1},
v))$ and the applications $\Pi_{f}$, $f\in\mC_{c}(\Sun\times\R{})$ by:
\[\begin{array}{cccl}
\Pi_{f}: & \Om & \tend & \bD([0, T], \R{})\\
& m &\mapsto & \cro{m}{f}.
\end{array}\]
Let $(P_{n})_{n}$ a sequence of probabilities on $\Om$ and $(\Pi_{f}P_{n})= P_{n}\circ \Pi_{f}^{-1} \in
\bD([0, T], \R{})$. We recall the following result:
\medskip
\begin{lem}
If for all $k\geq 0$, the sequence $(\Pi_{f_{k}}P_{n})_{n}$ is tight in $\mM_{1}(\bD([0,T], \R{}))$, then the
sequence $(P_{n})_{n}$ is tight in $\mM_{1}(\bD([0, T], (\mM_{1}, v)))$.
\end{lem}

Hence, it suffices to have a criterion for tightness in $\bD([0, T], \R{})$.
Let $X_{t}^{n}$ be a sequence of processes in $\bD([0, T], \R{})$ and $\mF_{t}^{n}$ a sequence of filtrations
such that $X^{n}$ is $\mF^{n}$-adapted. Let $\phi^{n}=\{\text{stopping times for }\mF^{n}\}$.
We have (cf. Billingsley \cite{cf:BillConvProb}):
\medskip
\begin{lem}[(Aldous' criterion)]\label{lem:aldous}If the following holds,\begin{enumerate}
	\item $\mL\left(\sup_{t\leq T}\left|X_{t}^{n}\right|\right)_{n}$ is tight,\label{it:mLtight}
	\item \label{it:SS'}For all $\eps>0$ and $\eta>0$, there exists $\delta>0$ such that\[\limsup_{n}
\sup_{S, S'\in \phi^{n}; S\leq S'\leq (S+\delta)\wedge T} \bP\left(\left|X_{S}^{n}-
X_{S'}^{n}\right|>\eta\right)\leq \eps,\]
\end{enumerate} then $\mL(X^{n})$ is tight.
\end{lem}
\medskip
\begin{prop}
The sequence $\mL(\nu^{N, (\om)})$ is tight in $\bD([0, T], (\mM_{F}, v))$.
\end{prop}

\begin{proof}
For all $\eps>0$, for all $k\geq 1$ (the case $k=0$ is straightforward),
	\[\bP\left(\sup_{t\leq T} \left|\cro{\nu^{N, (\om)}_{t}}{f_{k}}\right|>\frac{1}{\eps}\right)\leq \eps
\N{f_{k}}_{\infty}\bE\left[\sup_{t\leq T} \left|\underbrace{\cro{\nu^{N, (\om)}_{t}}{1}}_{= 1}\right|\right],
\quad \text{(Markov Inequality)}.\]
The tightness follows.

For all $k\geq 1$, we have the following decomposition:
\[\cro{\nu^{N, (\om)}_{t}}{f_{k}} = \cro{\nu^{N, (\om)}_{0}}{f_{k}} + A^{N, (\om)}_{t}(f_{k}) + M^{N,
(\om)}_{t}(f_{k}),\]
where $A^{N, (\om)}_{t}(f_{k})$ is a process of bounded variations, and $M^{N, (\om)}_{t}(f_{k})$ is a
square-integrable martingale. Then it suffices to verify Lemma \ref{lem:aldous}, (\ref{it:SS'}) for $A$ and
$M$ separately.
For all $\eps>0$ and $\eta>0$, for all stopping times $S, S'\in \phi^{N}; S\leq S'\leq (S+\delta)\wedge T$, we
have:
\begin{align*}
a_N&:=\bP\left(\left|A^{N, (\om)}_{S'}(f_{k}) - A^{N, (\om)}_{S}(f_{k})\right|>\eta\right),\\
&\leq \frac{1}{\eta} \bE\left[\int_{S}^{S'}{\d s\left|\cro{\nu^{N, (\om)}_{s}}{f_{k}'\cdot(b[\cdot, \nu^{N,
(\om)}_{s}] + c)}\right|}\right]+ \frac{1}{\eta}\bE\left[\frac{1}{2}\int_{S}^{S'}{\d s\left|\cro{\nu^{N,
(\om)}_{s}}{f_{k}''}\right|}\right],\\
&\leq\frac{C}{\eta} \bE\left[S'-S\right] \leq \eps,\quad\text{for $\delta$ sufficiently small.}
\end{align*}(we use here that $f_{k}$ are of compact support for $k\geq 1$; in particular the function
$(x,\pi)\mapsto f'_{k}(x,\pi)c(x, \pi)$ is bounded).
Furthermore,\begin{align*}
\bP\left(\left|M^{N, (\om)}_{S'}(f_{k}) - M^{N, (\om)}_{S}(f_{k})\right|>\eta\right) &= \bP\left(\left|M^{N,
(\om)}_{S'}(f_{k}) - M^{N, (\om)}_{S}(f_{k})\right|^{2}>\eta^{2}\right),\\
&\leq \frac{1}{\eta^{2}} \bE\left[\left|M^{N, (\om)}_{S'}(f_{k}) - M^{N,
(\om)}_{S}(f_{k})\right|^{2}\right],\\
&\leq \frac{1}{(N\eta)^{2}}\bE\left[\Som{i}{1}{N}{\int_{S}^{S'}{f'^{2}_{k}(x^{i}_{s}, \om_{i})\d s}}\right]
\leq \frac{C}{N\eta^{2}}\delta.\qedhere
\end{align*}
\end{proof}
At this point, $\mL(\nu^{N, (\om)})$ is tight in $\bD([0, T], (\mM_{F}, v))$.

\subsubsection{Equation satisfied by any accumulation point in $\bD([0, T], (\mM_{F}, v))$}
Using hypothesis \eqref{eq:H3nu}, it is easy to show that the following equation is satisfied for every
accumulation point $\nu$, for every $f\in\mC_{c}^{2}(\Sun\times\R{})$ (we use here that $\Sun$ is compact):
 \begin{align}
\cro{\nu_{t}}{f} &= \cro{\nu_{0}}{f} + \int_{0}^{t}{\d s\cro{\nu_{s}}{f'\cdot (b[\cdot, \nu_{s}]+c)}} +
\frac{1}{2}\int_{0}^{t}{\d s\cro{\nu_{s}}{f''}}.\label{eq:limnuom}
\end{align}

\medskip
For any accumulation point $\nu$, the following lemma gives a uniform bound for the second marginal of $\nu$:
\medskip
\begin{lem}\label{lem:marge}
Let $Q$ be an accumulation point of $\mL(\nu^{N, (\om)})_{N}$ in $\bD([0, T], (\mM_{1}, v))$ and let be
$\nu\sim Q$.
For all $t\in[0, T]$, we define by $(\nu_{t, 2})$ the second marginal of $\nu_{t}$:
\[\forall A\in\mB(\R{}),\quad (\nu_{t, 2})(A) = \int_{\Sun\times A}{\nu_{t}(\d x, \d\pi)}.\]
Then, for all $t\in[0,T]$,
\begin{equation*}
\label{eq:margnu}
\int_{\R{}}{\sup_{x\in \Sun}\left|c(x, \pi)\right|(\nu_{t, 2})(\d\pi)} \leq
\int_{\R{}}{\sup_{x\in\Sun}\left|c(x, \pi)\right|\mu(\d\pi)}.
\end{equation*}
\end{lem}

\begin{proof}
Let $\phi$ be a $\mC^{2}$ positive function such that $\phi\equiv 1$ on $[-1, 1]$, $\phi\equiv 0$ on $[-2, 2]$
and $\N{\phi}_{\infty}\leq 1$. Let, 
\[\forall k\geq 1,\quad \phi_{k}:=\pi\mapsto \phi\left(\frac{\pi}{k}\right).\]
Then $\phi_{k}\in\Cdc$ and $\phi_{k}(\pi)\tend_{k\tend\infty} 1$, for all $\pi$. We have also for all
$\pi\in\R{}$, $|\phi_{k}(\pi)|\leq \N{\phi}_{\infty}$, $|\phi'_{k}(\pi)|\leq \N{\phi'}_{\infty}$,
$|\phi''_{k}(\pi)|\leq \N{\phi''}_{\infty}$.

We have successively, denoting $S(\pi):= \sup_{x\in\Sun}\left|c(x, \pi)\right|$,
\begin{align}
\int_{\Sun\times \R{}}{S(\pi)\nu_{t}(\d x,\d\pi)} &= \int_{\Sun\times
\R{}}{\liminf_{k\tend\infty}\phi_{k}(\pi)S(\pi)\nu_{t}(\d x,\d\pi)},\nonumber\\
	&\leq \liminf_{k\tend\infty} \int_{\Sun\times \R{}}{\phi_{k}(\pi)S(\pi)\nu_{t}(\d x,\d\pi)},\
\text{(Fatou's lemma)},\nonumber\\
	&= \liminf_{k\tend\infty} \lim_{N\tend\infty}\int_{\Sun\times \R{}}{\phi_{k}(\pi)S(\pi)\nu^{N,
(\om)}_{t}(\d x,\d\pi)},\label{eq:ineqfatou}\\
	&\leq \lim_{N\tend\infty}\int_{\Sun\times \R{}}{S(\pi) \nu^{N, (\om)}_{t}(\d x,\d\pi)}, \text{(since
$\N{\phi}_{\infty}\leq 1$)}.\nonumber
\end{align}
The equality \eqref{eq:ineqfatou} is true since $(x, \pi) \mapsto \phi_{k}(\pi)S(\pi)$ is of compact support
in $\Sun\times\R{}$ (recall that $S$ is supposed to be continuous by hypothesis).

But, by definition of $\nu^{N, (\om)}_{t}$, and using the hypothesis \eqref{eq:Hmu} concerning $\mu$, we have,
\begin{equation}
\label{eq:lfgnom}
\lim_{N\tend\infty}\int_{\Sun\times \R{}}{S(\pi) \nu^{N, (\om)}_{t}(\d x,\d\pi)} = \int_{\R{}}{S(\pi)
\mu(\d\pi)}.
\end{equation}The result follows.\qedhere
\end{proof}

\begin{rem}
\eqref{eq:lfgnom} is only true for $\bbP$-almost every sequence $(\om)$. We assume that the sequence $(\om)$
given at the beginning satisfies this property.
\end{rem}
\subsubsection{Tightness in $\bD([0, T], (\mM_{F}, w))$}
We have the following lemma (cf. \cite{cf:RC}):
\medskip
\begin{lem}\label{theo:meleard}
Let $(X_{n})$ be a sequence of processes in $\bD([0, T], (\mM_{F}, w))$ and $X$ a process belonging to
$\mC([0, T], (\mM_{F}, w))$. Then,
\[X_{n}\stackrel{\mL}{\rightarrow} X \Leftrightarrow \left\{\begin{array}{ll} X_{n}\stackrel{\mL}{\rightarrow}
X &in\ \bD([0, T], (\mM_{F}, v)),\\ \cro{X_{n}}{1}\stackrel{\mL}{\rightarrow} \cro{X}{1} & in\ \bD([0, T],
\R{}).\end{array}\right.\]
\end{lem}

So, it suffices to show, for any accumulation point $\nu$:
\begin{enumerate}
	\item $\cro{\nu}{1}=1$: Eq. \eqref{eq:limnuom} is true for all $f\in\mC^{2}_{c}(\Sun\times\R{})$, so
in particular for $f_{k}(x, \pi):= \phi_{k}(\pi)$. Using the boundedness shown in lemma \ref{lem:marge}, we
can apply dominated convergence theorem in Eq. \eqref{eq:limnuom}. We then have $\cro{\nu_{t}}{1}=1$, for all
$t\in[0,T]$. The fact that Eq. \eqref{eq:limnuom} is verified for all $f\in \mC^{2}_{b}(\Sun\times\R{})$ can
be shown in the same way.
	\item Continuity of the limit: For all $0\leq s\leq t\leq T$, for all
$f\in\mC^{2}_{b}(\Sun\times\R{})$,
\begin{align*}
	\left|\cro{\nu_{t}}{f}-\cro{\nu_{s}}{f}\right|&\leq K\int_{s}^{t}{\left|\cro{\nu_{u}}{f'\cdot b[\cdot
, \nu_{u}]}\right|\d u} + \frac{1}{2}\int_{s}^{t}{\left|\cro{\nu_{u}}{f''}\right|\d u}\\&+
	\int_{s}^{t}{\left|\cro{\nu_{u}}{f'\cdot c}\right|\d u} \leq C\times \left|t-s\right|, \quad\text{for
some constant $C$.}
\end{align*} Noticing that we used again Lemma \ref{lem:marge} to bound the last term, we have the result.
\end{enumerate}

\begin{rem}It is then easy to see that the second marginal (on the disorder) of any accumulation point is
$\mu$.\end{rem}

At this point $\mL(\nu^{N, (\om)})$ is tight in $\bD([0, T], (\mM_{F}, w))$.
It remains to show the uniqueness of any accumulation point: it shows firstly that the sequence effectively
converges and that the limit does not depend on the given sequence $(\om)$.
\subsection{Uniqueness of the limit}
\medskip
\begin{prop}
\label{prop:uniciteMKV}
There exists a unique element $P$ of $\bD([0, T], \mM_{F}(\Sun\times\R{}))$ which satisfies equation
\eqref{eq:limnuom}, $P_{0}\in \mM_{1}$ and $P_{0, 2} =\mu$.
\end{prop}
\medskip
\begin{rem}
Since this proof is an adaptation of Oelschläger \cite{cf:Oelsch} Lemma 10, p.474, we only sketch the proof:
\end{rem}
We can rewrite Eq. \eqref{eq:limnuom} in a more compact way:
\begin{equation}
\label{eq:uninuLd}
\cro{P_{t}}{f} = \cro{P_{0}}{f} + \int_{0}^{t}{\cro{P_{s}}{L(P_{s})(f)}\d s}, \quad \forall
f\in\mC^{2}_{b}(\Sun\times\R{}), 0\leq t\leq T,
\end{equation}where,
\begin{equation*}
\label{eq:Lnud}
L(P)(f)(x,\om):= \frac{1}{2} f''(x,\om) + h(x, \om, P) f'(x,\om),\quad \forall x\in\Sun, \forall\om\in\R{},
\end{equation*}
and,
\begin{equation*}
\label{eq:bnud}
h(x, \pi, P) = b[x, P] + c(x, \pi).
\end{equation*}

Let $t\mapsto P_{t}$ be any solution of Eq. \eqref{eq:uninuLd}. One can then introduce the following SDE
(where $\xi\in\R{}$ and $\txi\in\Sun$ is its projection on $\Sun$):

\begin{equation}
\label{eq:sdeunid}
\left\{\begin{array}{ccl}
\d\xi_{t} &=& h(\txi_{t}, \om_{t}, P_{t}) \d t + \d W_{t}, \quad \xi_{0} = \txi_{0},\ (\txi_{0}, \om_{0}) \sim
P_{0},\\
\d\om_{t} &=& 0.
\end{array}\right.
\end{equation}
Eq. \eqref{eq:sdeunid} has a unique (strong) solution $(\xi_{t}, \om_{t})_{t\in[0, T]}=(\xi_{t},
\om_{0})_{t\in[0, T]}$. The proof of uniqueness in Eq. \eqref{eq:limnuom} consists in two steps:

\begin{enumerate}
	\item $P_{t} = \mL(\txi_{t},\om_{t})$, for all $t\in[0, T]$,
	\item Uniqueness of $\txi$.
\end{enumerate}

We refer to Oelschläger \cite{cf:Oelsch} for the details. Another proof of uniqueness can be found in
\cite{cf:dPdH} or in \cite{cf:Gartner} via a martingale argument.

\section{Proof of the fluctuations result}
\label{sec:fluct}
Now we turn to the proof of Theorem \ref{theo:fluctquenchedinit}. To that purpose, we need to introduce some
distribution spaces:

\subsection{Distribution spaces}

Let $\mS:=\mS(\Sun\times\R{})$ be the usual Schwartz space of rapidly decreasing infinitely differentiable
functions. Let $\mD_{p}$ be the set of all differential operators of the form
$\partial_{u^{k}}\partial_{\pi^{l}}$ with $k+l\leq p$. We know from Gelfand and Vilenkin \cite{cf:GelfVil} p.
82-84, that we can introduce on $\mS$ a nuclear Fr\'echet topology by the system of seminorms $\N{\cdot}_{p}$,
$p=1,2,\dots,$ defined by
\[\N{\phi}_{p}^{2} = \Som{k}{0}{p}{\int_{\Sun\times\R{}}(1+|\pi|^{2})^{2p} \sum_{D\in\mD_{k}}{|D\phi(y,
\pi)|^{2}}\d y\d\pi}.\]

Let $\mS'$ be the corresponding dual space of tempered distributions. Although, for the sake of simplicity, we
will mainly consider $\eta^{N, (\om)}$ as a process in $\mC([0, T], \mS')$, we need some more precise
estimations to prove tightness and convergence. We need here the following norms:

For every integer $j$, $\alpha\in\R{+}$, we consider the space of all real functions $\vphi$ defined on
$\Sun\times\R{}$ with derivative up to order $j$ such that \[\N{\vphi}_{j,\alpha}:=
\left(\sum_{k_{1}+k_{2}\leq j}{\int_{\Sun\times\R{}}{\frac{|\partial_{x^{k_{1}}}\partial_{\pi^{k_{2}}}\vphi(x,
\pi)|^{2}}{1+|\pi|^{2\alpha}}\d x \d\pi}}\right)^{1/2} <\infty.\]
 Let $W^{j, \alpha}_{0}$ be the completion of $\mC_{c}^{\infty}(\Sun\times\R{})$ for this norm; $(W^{j,
\alpha}_{0}, \N{\cdot}_{j, \alpha})$ is a Hilbert space. Let $W_{0}^{-j, \alpha}$ be its dual space.

Let $C^{j,\alpha}$ be the space of functions $\vphi$ with continuous partial derivatives up to order $j$ such
that \[\lim_{|\pi|\tend\infty} \sup_{x\in\Sun} \frac{|\partial_{x^{k_{1}}}\partial_{\pi^{k_{2}}}\vphi(x,
\pi)|}{1+|\pi|^{\alpha}} = 0, \ \text{for all $k_{1}+k_{2}\leq j$,}\]
with norm \[\N{\vphi}_{C^{j,\alpha}} = \sum_{k_{1}+k_{2}\leq j} \sup_{x\in\Sun}
\sup_{\pi\in\R{}}\frac{|\partial_{x^{k_{1}}}\partial_{\pi^{k_{2}}}\vphi(x, \pi)|}{1+|\pi|^{\alpha}}.\]
We have the following embeddings:\[W_{0}^{m+j, \alpha} \hookrightarrow C^{j, \alpha}, m>1, j\geq 0, \alpha\geq
0,\] i.e. there exists some constant $C$ such that\begin{equation}
\label{eq:injsobo1}
\N{\vphi}_{C^{j, \alpha}}\leq C \N{\vphi}_{m+j, \alpha}.
\end{equation}
Moreover,
\[W_{0}^{m+j,\alpha}\hookrightarrow W_{0}^{j, \alpha+ \beta}, m>1, j\geq 0, \alpha\geq 0, \beta >1.\] Thus
there exists some constant $C$ such that 
\[\N{\vphi}_{j,\alpha+\beta} \leq C \N{\vphi}_{m+j, \alpha}.\]
We then have the following dual continuous embedding:
\begin{equation}
\label{eq:injsobo3}
W_{0}^{-j, \alpha+\beta} \hookrightarrow W_{0}^{-(m+j), \alpha}, \ m>1, \alpha\geq 0, \beta>1.
\end{equation}

It is quite clear that $\mS \hookrightarrow W_{0}^{j, \alpha}$ for any $j$ and $\alpha$, with a continuous
injection.

We now prove some continuity of linear mappings in the corresponding spaces:
\begin{lem}
For every $x, y \in\Sun$, $\om\in\R{}$, for all $\alpha$, the linear mappings $W_{0}^{3, \alpha} \tend \R{}$
defined by \[D_{x,y,\om}(\vphi):= \vphi(x,\om) -\vphi(y,\om); D_{x, \om}:= \vphi(x, \om); H_{x, \om} =
\vphi'(x, \om),\]are continuous and
\begin{align}
	\N{D_{x,y,\om}}_{-3, \alpha} \leq C|x-y| \left(1+ |\om|^{\alpha}\right),\label{eq:NDxyo}\\
	\N{D_{x,\om}}_{-3, \alpha} \leq C\left(1+ |\om|^{\alpha}\right),\label{eq:NDxo}\\
	\N{H_{x,\om}}_{-3, \alpha} \leq C\left(1+ |\om|^{\alpha}\right).\label{eq:NHxo}
\end{align}
\end{lem}

\begin{proof}
Let $\vphi$ be a function of class $\mC^{\infty}$ with compact support on $\Sun\times\R{}$, then,
\begin{align*}
	|\vphi(x, \om) - \vphi(y, \om)| &\leq |x-y| \sup_{u} \left|\vphi'(u, \om)\right|,\\
	&\leq |x-y| \left(1+ |\om|^{\alpha}\right) \sup_{u,\om}\left(\frac{\left|\vphi'(u,
\om)\right|}{1+|\om|^{\alpha}}\right),\\
	&\leq |x-y| \left(1+ |\om|^{\alpha}\right) \N{\vphi}_{1,\alpha},\\
	&\leq C|x-y| \left(1+ |\om|^{\alpha}\right) \N{\vphi}_{3, \alpha},
\end{align*}
following \eqref{eq:injsobo1} with $j=1$ and $m=2>1$. Then, \eqref{eq:NDxyo} follows from a density argument.
\eqref{eq:NDxo} and \eqref{eq:NHxo} are proved in the same way.
\end{proof}

\subsection{The non-linear process} The proof of convergence is based on the existence of the non-linear
process associated to McKean-Vlasov equation. Such existence has been studied by numerous authors (eg. Dawson
\cite{Dawson1983}, Jourdain-M\'el\'eard \cite{Jourdain1998}, Malrieu \cite{Malrieu2003}, Shiga-Tanaka
\cite{Shiga1985}, Sznitman \cite{cf:Sznit84}, \cite{cf:SznitSflour}) mostly in order to prove some propagation
of chaos properties in systems without disorder. We consider the present similar case where disorder is
present. Let us give some intuition of this process. One can replace the non-linearity in Eq. \eqref{eq:MKV1L}
by an arbitrary measure $m(\d x, \d\om)$:
\begin{equation*}
\partial_{t} q_{t}^{\om} = \frac{1}{2} \partial_{xx}q_{t}^{\om} - \partial_{x}\left[\left(b[x,m]+ c(x,
\om)\right)q_{t}^{\om}\right].
\label{eq:MKV1m}
\end{equation*}
In this particular case, it is usual to interpret $q_{t}^{\om}$ as the time marginals of the following
diffusion:
\[dx_{t}^{\om} = dB_{t} + b[x_{t}^{\om}, m] dt + c(x_{t}^{\om}, \om) dt,\ \om\sim\mu.\] 
It is then natural to consider the following problem, where $m$ is replaced by the proper measure $P$: on a
filtered probability space $(\Om, \mF, \mF_{t}, B, x_{0}, Q)$, endowed with a Brownian motion $B$ and with a
$\mF_{0}$ measurable random variable $x_{0}$, we introduce the following system:
\begin{equation}
\label{eq:pnonlin}
\left\{\begin{array}{ccl}
 x_{t}^{\om} &=& x_{0} + \int_{0}^{t}{b[x_{s}^{\om}, P_{s}]\d s} + \int_{0}^{t}{c(x_{s}^{\om}, \om)\d s} +
B_{t},\\
 \om&\sim& \mu,\\
 P_{t} &=& \mL(x_{t}, \om), \forall t\in[0, T].
 \end{array}\right.\end{equation}
\begin{prop}
There is pathwise existence and uniqueness for Equation \eqref{eq:pnonlin}.
\end{prop}
\begin{proof}
The proof is the same as given in Sznitman \cite{cf:SznitSflour}, Th 1.1, p.172, up to minor modifications.
The main idea consists in using a Picard iteration in the space of probabilities on $
\mC([0,T],\Sun\times\R{})$ endowed with an appropriate Wasserstein metric. We refer to it for details.
\end{proof}
\subsection{Fluctuations in the quenched model}
The key argument of the proof is to explicit the speed of convergence as $N\tend\infty$ for the rotators to
the non-linear process (see Prop. \ref{prop:estimpnonlin}).

A major difference between this work and \cite{cf:FernMeleard} is that, since in our quenched model, we only
integrate w.r.t. oscillators \emph{and not} w.r.t. the disorder, one has to deal with remaining terms, see
$Z_{N}$ in Proposition \ref{prop:estimpnonlin}, to compare with \cite{cf:FernMeleard}, Lemma 3.2, that would
have disappeared in the \emph{averaged model}. The main technical difficulty of Proposition
\ref{prop:estimpnonlin} is to control the asymptotic behaviour of such terms, see \eqref{eq:estimZN}. As in
\cite{cf:FernMeleard}, having proved Prop. \ref{prop:estimpnonlin}, the key argument of the proof is a uniform
estimation of the norm of the process $\eta^{N, (\om)}$, see Propositions \ref{prop:estimeta1} and
\ref{prop:estimsupetat}, based on the generalized stochastic differential equation verified by $\eta^{N
(\om)}$, see \eqref{eq:itoetaN}.
\subsubsection{Preliminary results}
\label{sec:premresfluct}
We consider here a fixed realization of the disorder $(\om) = (\om_{1}, \om_{2}, \dots)$. On a common filtered
probability space $(\Om, \mF, \mF_{t}, (B^{i})_{i\geq1}, Q)$, endowed with a sequence of i.i.d.
$\mF_{t}$-adapted Brownian motions $(B_{i})$ and with a sequence of i.i.d. $\mF_{0}$ measurable random
variables $(\xi^{i})$ with law $\lambda$, we define as $x^{i, N}$ the solution of \eqref{eq:diffbc}, and as
$x^{\om_{i}}$ the solution of \eqref{eq:pnonlin}, with the same Brownian motion $B^{i}$ and with the same
initial value $\xi^{i}$.

The main technical proposition, from which every norm estimation of $\eta^{N, (\om)}$ follows is the
following:
\begin{prop}
\label{prop:estimpnonlin}
\begin{equation}
\label{eq:espxin}
\bE\left[\sup_{t\leq T}\left|x_{t}^{i,N} - x_{t}^{\om_{i}}\right|^{2}\right] \leq C/N + Z_{N}(\om_{1},\dots,
\om_{N}),
\end{equation}
where the random variable $(\om)\mapsto Z_{N}(\om)$ is such that:
\begin{equation}
\label{eq:estimZN}
\lim_{A\tend\infty} \limsup_{N\tend\infty}\bbP\left(N Z_{N}(\om)>A\right) = 0.
\end{equation}
\end{prop}
The (rather technical) proof of Proposition \ref{prop:estimpnonlin} is postponed to the end of the document
(see \S \ref{app:propZN}).
Once again, we stress the fact that the term $Z_N$ would have disappeared in the averaged model.

The first norm estimation of the process $\eta^{N, (\om)}$ (which will be used to prove tightness) is a direct
consequence of Proposition \ref{prop:estimpnonlin} and of a Hilbertian argument:
\medskip
\begin{prop}
\label{prop:estimeta1}
Under the hypothesis \eqref{eq:Hmufluct} on $\mu$, the process $\eta^{N, (\om)}$ satisfies the following
property: for all $T>0$,
\begin{equation}
\label{eq:estimeta1}
\sup_{t\leq T} \bE\left[\N{\eta_{t}^{N, (\om)}}^{2}_{-3, 2\alpha}\right] \leq A_{N}(\om_{1}, \dots, \om_{N}),
\end{equation}
where \begin{equation*}
\label{eq:estimAN}
\lim_{A\tend\infty}\limsup_{N\tend\infty}\bbP\left(A_{N}>A\right)=0.
\end{equation*}
\end{prop}

\begin{proof}
For all $\vphi\in W_{0}^{3, 2\alpha}$, writing 
\begin{align*}
\cro{\eta_{t}^{N, (\om)}}{\vphi} &= \frac{1}{\sqrt{N}}\Som{i}{1}{N}{\left\{\vphi(x_{t}^{i,
N},\om_{i})-\vphi(x_{t}^{\om_{i}},\om_{i})\right\}}+
\frac{1}{\sqrt{N}}\Som{i}{1}{N}{\left\{\vphi(x_{t}^{\om_{i}},\om_{i})-\cro{P_{s}}{\vphi}\right\}},\\ 
&=:S_{t}^{N, (\om)}(\vphi) + T_{t}^{N, (\om)}(\vphi),
\end{align*}
we have:\begin{equation}
\label{eq:etaST}
\cro{\eta_{t}^{N, (\om)}}{\vphi}^{2} \leq 2 \left(S_{t}^{N, (\om)}(\vphi)^{2} + T_{t}^{N,
(\om)}(\vphi)^{2}\right).
\end{equation}
But, by convexity, \[S_{t}^{N, (\om)}(\vphi)^{2}\leq \Som{i}{1}{N}{D^{2}_{x_{t}^{i, N}, x_{t}^{\om_{i}},
\om_{i}}(\vphi)}.\] Then, applying the latter equation to an orthonormal system $(\vphi_{p})_{p\geq 1}$ in the
Hilbert space $W_{0}^{3, 2\alpha}$, summing on $p$, we have by Parseval's identity on the continuous
functional $D_{x_{t}^{i, N}, x_{t}^{\om_{i}}}$,
\begin{align}
	\bE\left[\N{S_{t}^{N, (\om)}}^{2}_{-3, 2\alpha}\right] &\leq \bE\left[\Som{i}{1}{N}{\N{D_{x_{t}^{i,
N}, x_{t}^{\om_{i}}, \om_{i}}}_{-3, 2\alpha}^{2}}\right],\nonumber\\
	&\leq C\Som{i}{1}{N}{(1+|\om_{i}|^{4\alpha})\bE\left[\left|x_{t}^{i, N} -
x_{t}^{\om_{i}}\right|^{2}\right]},\label{eq:etaSD}\\
	&\leq C \Som{i}{1}{N}{\left(1+ |\om_{i}|^{4\alpha}\right)} \left(C/N + Z_{N}(\om_{1}, \dots,
\om_{N})\right),\label{eq:etaSZ}
\end{align}where we used \eqref{eq:NDxyo} in \eqref{eq:etaSD}, and \eqref{eq:espxin} in \eqref{eq:etaSZ}.

On the other hand,
\begin{align*}
	\bE\left[T_{t}^{N, (\om)}(\vphi)^{2}\right] &=  \frac{1}{N}
\bE\left[\left\{\Som{i}{1}{N}{(\vphi(x_{t}^{\om_{i}}, \om_{i}) - \cro{P_{t}}{\vphi})}\right\}^{2}\right],\\
	&= \frac{1}{N} \bE\left[\Som{i}{1}{N}{(\vphi(x_{t}^{\om_{i}}, \om_{i}) -
\cro{P_{t}}{\vphi})^{2}}\right]\\
	&+ \frac{1}{N} \bE\left[\sum_{i\neq j}(\vphi(x_{t}^{\om_{i}}, \om_{i}) -
\cro{P_{t}}{\vphi})(\vphi(x_{t}^{\om_{j}}, \om_{j}) - \cro{P_{t}}{\vphi})\right],\\
	&\leq \frac{2}{N} \bE\left[\Som{i}{1}{N}{(\vphi(x_{t}^{\om_{i}}, \om_{i})^{2} +
\cro{P_{t}}{\vphi}^{2})}\right] + \frac{1}{N} \sum_{i\neq j}G(\vphi)(\om_{i})G(\vphi)(\om_{j}),\\
	&\leq \frac{2}{N} \bE\left[\Som{i}{1}{N}{\vphi(x_{t}^{\om_{i}}, \om_{i})^{2}}\right] + 
2\cro{P_{t}}{\vphi}^{2} + \left(\frac{1}{\sqrt{N}}\Som{i}{1}{N}{G(\vphi)(\om_{i})}\right)^{2},
\end{align*}
where $G(\vphi)(\om):=\int\vphi(y, \om_{i})P_{t}^{\om_{i}}(\d y) - \cro{P_{t}}{\vphi}$. If we apply the same
Hilbertian argument as for $S^{N, (\om)}$, we see
\begin{equation}
\label{eq:estimTNfinal}
\bE\left[\N{T_{t}^{N, (\om)}}_{-3, 2\alpha}^{2}\right] \leq \frac{2C}{N}
\bE\left[\Som{i}{1}{N}{(1+|\om_{i}|^{4\alpha})}\right] +  C +
\N{\phi\mapsto\left(\frac{1}{\sqrt{N}}\Som{i}{1}{N}{G(\phi)(\om_{i})}\right)}_{-3, 2\alpha}^{2},
\end{equation}It is easy to see that the last term in \eqref{eq:estimTNfinal} can be reformulated as
$B_{N}(\om_{1}, \dots, \om_{N})$, with the property that
$\lim_{A\tend\infty}\limsup_{N\tend\infty}\bbP(B_{N}>A)=0$. Combining \eqref{eq:estimZN}, \eqref{eq:etaST}, 
\eqref{eq:etaSZ} and \eqref{eq:estimTNfinal}, Proposition \ref{prop:estimeta1} is proved.
\end{proof}
\subsubsection{Tightness of the fluctuations process}

Applying Ito's formula to \eqref{eq:diffbc}, we obtain, for all $\vphi$ bounded function on $\Sun\times\R{}$,
with two bounded derivatives w.r.t. $x$, for every sequence $(\om)$, for all $t\leq T$:
\begin{equation}
\label{eq:itoetaN}
\cro{\eta_{t}^{N, (\om)}}{\vphi} = \cro{\eta_{0}^{N, (\om)}}{\vphi} + \int_{0}^{t}{\cro{\eta_{s}^{N,
(\om)}}{\mL_{s}^{\nu^{N}}(\vphi)}\d s} + M_{t}^{N, (\om)}(\vphi),
\end{equation} where, for all $y\in\Sun$, $\pi\in\R{}$,
\[\mL_{s}^{\nu^{N}}(\vphi)(y, \pi) = \frac{1}{2} \vphi''(y, \pi) + \vphi'(y, \pi)\left(b[y, \nu_{s}^{N}] +
c(y, \pi)\right) + \cro{P_{s}}{\vphi'(\cdot, \cdot) b(\cdot, y, \pi)},\] and $M_{t}^{N, (\om)}(\vphi)$ is a
real continuous martingale with quadratic variation process
\[\left\langle M^{N, (\om)}(\vphi)\right\rangle_{t} = \int_{0}^{t}{\cro{\nu_{s}^{N, (\om)}}{\vphi'(y,
\pi)^{2}}\d s}.\]

\begin{lem}
\label{lem:estimLmu}
For every $N$, the operator $\mL_{s}^{\nu^{N}}$ defines a linear mapping from $\mS$ into $\mS$ and for all
$\vphi\in \mS$,
\[\N{\mL_{s}^{\nu^{N}}(\vphi)}_{3, 2\alpha}^{2} \leq C \N{\vphi}_{6,\alpha}^{2}.\]
\end{lem}

\begin{proof}
The terms $\frac{1}{2}\vphi''(y,\pi)$ and $\vphi'(y,\pi) b[y, \nu_{s}^{N}]$ clearly satisfy the lemma. We
study the two remaining terms:
\begin{align*}
	\N{\cro{P_{s}}{\vphi'b(\cdot,y,\pi)}}_{3, 2\alpha}^{2} &= \sum_{k_{1}+k_{2}\leq
3}\int_{\Sun\times\R{}}{
\frac{\cro{P_{s}}{\vphi'\partial_{y^{k_{1}}}\partial_{\pi^{k_{2}}}b(\cdot,y,\pi)}^{2}}{1+ |\pi|^{4\alpha}}\d
y\d\pi},\\
	&\leq  C \int_{\R{}}{\frac{1}{1+|\pi|^{4\alpha}}\d\pi} \int_{\Sun\times\R{}}{\vphi'(y,
\pi)^{2}P_{s}(\d y,\d\pi)},\\
	&\leq C \N{\vphi}_{C^{3, \alpha}}^{2}\int_{\R{}}{\frac{1}{1+|\pi|^{4\alpha}}\d\pi}
\int_{\Sun\times\R{}}{(1+ |\pi|^{\alpha})^{2}P_{s}(\d y,\d\pi)},\\
	&\leq C \N{\vphi}_{6, \alpha}^{2}\int_{\R{}}{\frac{1}{1+|\pi|^{4\alpha}}\d\pi} \int_{\R{}}{(1+
|\pi|^{\alpha})^{2}\mu(\d\pi)}.
\end{align*}
And,\begin{align*}
	\N{\vphi'(y, \pi)c(y,\pi)}_{3, 2\alpha}^{2} &= \sum_{k_{1}+k_{2}\leq 3}\int_{\Sun\times\R{}}{
\frac{\left(\partial_{y^{k_{1}}}\partial_{\pi^{k_{2}}}\left\{\vphi'(y, \pi)c(y,\pi)\right\}\right)^{2}}{1+
|\pi|^{4\alpha}}\d y\d\pi}.\\
\end{align*}
It suffices to estimate, for every differential operator $D_{i} = \partial_{y^{u_{i}}}\partial_{\pi^{v_{i}}}$,
$i=1,2$ with $u_{1}+u_{2}+v_{1}+v_{2}\leq3$, the following term:
\begin{align*}
	\int_{\Sun\times \R{}}{\frac{|D_{1}\vphi'(y,\pi)D_{2}c(y, \pi)|^{2}}{1+|\pi|^{4\alpha}}\d y\d\pi}&\leq
\int_{\Sun\times \R{}}\frac{|D_{1}\vphi'(y,\pi)|^{2}}{(1+|\pi|^{\alpha})^{2}}\frac{|D_{2}c(y,
\pi)|^{2}(1+|\pi|^{\alpha})^{2}}{1+|\pi|^{4\alpha}}\d y\d\pi,\\
	&\leq C \N{\vphi}_{6, \alpha}^{2} \int_{\R{}}{\frac{\sup_{y\in\Sun}|D_{2}c(y,
\pi)|^{2}}{1+|\pi|^{2\alpha}}\d\pi}.
\end{align*}
The result follows from the assumptions made on $c$.
\end{proof}

For the tightness criterion used below, we need to ensure that the trajectories of the fluctuations process
are almost surely continuous: in that purpose, we need some more precise evaluations than in Prop.
\ref{prop:estimeta1}.
\begin{prop}
\label{prop:estimsupMt}
The process $(M_{t}^{N, (\om)})$ satisfies, for every $(\om)$, and for every $T>0$,
\begin{equation*}
\label{eq:estimsupMt}
\bE\left[\sup_{t\leq T} \N{M_{t}^{N, (\om)}}^{2}_{-3, 2\alpha}\right] \leq \frac{C}{N} \Som{i}{1}{N}{\left(1+
|\om_{i}|^{4\alpha}\right)}.
\end{equation*}
\end{prop}
\begin{rem}
\label{rem:remsupNMt}
In particular, a consequence of \eqref{eq:Hmufluct} is that, for $\bbP$-almost every sequence
$(\om)$,\begin{equation}
\label{eq:remsupNMt}
\sup_{N}\bE\left[\sup_{t\leq T} \N{M_{t}^{N, (\om)}}^{2}_{-3, 2\alpha}\right] \leq \sup_{N}\frac{C}{N}
\Som{i}{1}{N}{\left(1+ |\om_{i}|^{4\alpha}\right)}<\infty.
\end{equation}
\end{rem}
\begin{proof}
Let $(\vphi_{p})_{p\geq 1}$ a complete orthonormal system in $W_{0}^{3, 2\alpha}$. For fixed $N$, by Doob's
inequality, $\sum_{p\geq1} \bE\left[\sup_{t\leq T}(M_{t}^{N, (\om)}(\vphi_{p}))^{2}\right]$ is bounded by
\begin{align*}
	C\sum_{p\geq 1} \bE\left[M_{T}^{N, (\om)}(\vphi_{p})^{2}\right] &= C \sum_{p\geq 1}
\bE\left[\int_{0}^{T}{\cro{\nu_{s}^{N, (\om)}}{\vphi_{p}'(y, \pi)^{2}}\d s}\right],\\
	&= \frac{1}{N} \Som{i}{1}{N}{\bE\left[\int_{0}^{T}{\sum_{p\geq 1} \vphi_{p}'(x_{s}^{i, N},
\om_{i})^{2}\d s}\right]},\\
	&= \frac{1}{N} \Som{i}{1}{N}{\bE\left[\int_{0}^{T}{\N{H_{x_{s}^{i, N}, \om_{i}}}_{3, 2\alpha}^{2}\d
s}\right]},\\
	&\leq \frac{C}{N} \Som{i}{1}{N}{\left(1+ |\om_{i}|^{4\alpha}\right)}, \text{(using
\eqref{eq:NHxo})}.\qedhere
\end{align*}
\end{proof}

\begin{prop}
\label{prop:estimsupetat}
For every $N$, every $(\om)$,
\begin{equation}
\label{eq:estimsupetat}
\bE\left[\sup_{t\leq T} \N{\eta_{t}^{N, (\om)}}^{2}_{-6, \alpha}\right] < C_{N}(\om_{1}, \dots, \om_{N}),
\end{equation}
with \begin{equation*}
\label{eq:estimsupetatCN}
\lim_{A\tend\infty}\limsup_{N\tend\infty}\bbP\left(C_{N}>A\right)=0.
\end{equation*}
\end{prop}
\begin{proof}
Let $(\psi_{p})$ be a complete orthonormal system in $W_{0}^{6, \alpha}$ of $\mC^{\infty}$ function on
$\Sun\times \R{}$ with compact support. We prove the stronger result:
\[\bE\left[\sum_{p\geq 1} \sup_{t\leq T} \cro{\eta_{t}^{N, (\om)}}{\psi_{p}}^{2}\right]<\infty.\]
Indeed, \[\cro{\eta_{t}^{N, (\om)}}{\psi_{p}}^{2}\leq C \left(\cro{\eta_{0}^{N, (\om)}}{\psi_{p}}^{2} + T
\int_{0}^{t}{\cro{\eta_{s}^{N, (\om)}}{\mL_{s}^{\nu^{N}}(\psi_{p})}^{2}\d s} + M^{N,
(\om)}_{t}(\psi_{p})^{2}\right).\] By Doob's inequality,
	\[\begin{split}\bE\left[\sum_{p\geq 1} \sup_{t\leq T} \cro{\eta_{t}^{N,
(\om)}}{\psi_{p}}^{2}\right]&\leq C\left(\bE\left[\N{\eta_{0}^{N, (\om)}}^{2}_{-6, \alpha}\right] +
\bE\int_{0}^{T}{\sum_{p\geq 1}\cro{\eta_{s}^{N, (\om)}}{\mL_{s}^{\nu^{N}}(\psi_{p})}^{2}\d s}\right.\\
	   &+ \left.\displaystyle \sum_{p\geq 1}\bE\left[M_{T}^{N,
(\om)}(\psi_{p})^{2}\right]\right).\end{split}\]
By Lemma \ref{lem:estimLmu}, we have:
\[\left|\cro{\eta_{s}^{N, (\om)}}{\mL_{s}^{\nu^{N}}(\psi)}\right|\leq C \N{\eta_{s}^{N, (\om)}}_{-3,
2\alpha}\N{\psi}_{6, \alpha}.\]
Then, 
\begin{align*}
	\bE \left[\int_{0}^{T}{\sum_{p\geq 1}\cro{\eta_{s}^{N, (\om)}}{\mL_{s}^{\nu^{N}}(\psi_{p})}^{2}\d
s}\right]&\leq  C^{2}\int_{0}^{T}{\bE\left[\N{\eta_{s}^{N, (\om)}}_{-3, 2\alpha}^{2}\right]\d s},\\
	&\leq C^{2}T \sup_{s\leq T} \bE\left[\N{\eta_{s}^{N, (\om)}}_{-3, 2\alpha}^{2}\right] \leq C^{2}T
A_{N},
\end{align*}
where $A_{N}$ is defined in Proposition \ref{prop:estimeta1}. The result follows.
\end{proof}
\begin{prop}
\label{prop:contWeta}
\begin{enumerate}
	\item For every $N$, for $\bbP$-almost every $(\om)$, the trajectories of the fluctuations process
$\eta^{N, (\om)}$ are almost surely continuous in $\mS'$,
	\item For every $N$, for $\bbP$-almost every $(\om)$, the trajectories of $M^{N, (\om)}$ are almost
surely continuous in $\mS'$.
\end{enumerate}
\end{prop}
\begin{proof}
	We only prove for $M^{N, (\om)}$, since, using Proposition \ref{prop:estimsupetat}, the proof is the
same for $\eta^{N, (\om)}$. Let $(\vphi_{p})$ be a complete orthonormal system in $W_{0}^{-3, 2\alpha}$, then
for every fixed $N$ and $(\om)$, we know from the proof of Proposition \ref{prop:estimsupMt}, that for all
$\eps>0$, there exists some $M_{0}>0$ such that \[\sum_{p\geq M_{0}} \sup_{t\leq T} (M_{t}^{N,
(\om)}(\vphi_{p}))^{2}<\frac{\eps}{3}, a.s.\] Let $(t_{m})$ be a sequence in $[0,T]$ such that
$t_{m}\tend_{m\tend\infty} t$.
\begin{align*}
	\N{M_{t_{m}}^{N,(\om)}- M_{t}^{N,(\om)}}_{-3, 2\alpha}^{2} &= \sum_{p\geq 1}
\left(M_{t_{m}}^{N,(\om)}- M_{t}^{N,(\om)}\right)^{2}(\vphi_{p}),\\
	&\leq \Som{p}{1}{M_{0}}{\left(M_{t_{m}}^{N,(\om)}- M_{t}^{N,(\om)}\right)^{2}(\vphi_{p})} +
\frac{2\eps}{3} \leq \eps,
\end{align*}
if $t_{m}$ is sufficiently large.\qedhere
\end{proof}

We are now in position to prove the tightness of the fluctuations process. Let us recall some notations: for
fixed $N$ and $(\om)$ $\mH^{N,(\om)}$ is the law of the process $\eta^{N, (\om)}$. Hence, $\mH^{N, (\om)}$ is
an element of $\mM_{1}(\mC([0, T], \mS'))$, endowed with the topology of weak convergence and with $\mB^{*}$,
the smallest $\sig$-algebra such that the evaluations $Q\mapsto \cro{Q}{f}$ are measurable, $f$ being
measurable and bounded.

We will denote by $\Theta_{N}$ the law of the random variable $(\om)\mapsto \mH^{N, (\om)}$.
The main result of this part is the following:

\begin{theo}
\label{theo:tightthetaxi}
\begin{enumerate}
	\item for $\bbP$-almost every sequence $(\om)$, the law of the process $M^{N, (\om)}$ is tight in
$\mM_{1}(\mC([0, T], \mS'))$,
	\item The law of the sequence $(\om)\mapsto \mH^{N,(\om)}$ is tight on $\mM_{1}\left(\mM_{1}(\mC([0,
T], \mS'))\right)$.
\end{enumerate}
\end{theo}

Before proving Theorem \ref{theo:tightthetaxi}, we recall the following result and notations (cf. Mitoma
\cite{cf:Mitomatight}, Th 3.1, p. 993):

\begin{prop}[(Mitoma's criterion)]
Let $(P_{N})$ be a sequence of probability measures on $(\mC_{\mS'}:=\mC([0, T], \mS'), \mB_{C_{\mS'}})$. For
each $\vphi\in \mS$, we denote by $\Pi_{\vphi}$ the mapping of $\mC_{\mS'}$ to $\mC:=\mC([0, T], \R{})$
defined by \[\Pi_{\vphi}: \psi(\cdot)\in\mC_{\mS'}\mapsto \cro{\psi(\cdot)}{\vphi} \in \mC.\]Then, if for all
$\vphi\in\mS$, the sequence $(P_{N}\Pi_{\vphi}^{-1})$ is tight in $\mC$, the sequence $(P_{N})$ is tight in
$\mC_{\mS'}$.
\end{prop}
\begin{rem}
\label{rem:phisep}
A closer look to the proof of Mitoma shows that it suffices to verify the tightness of
$(P_{N}\Pi_{\vphi}^{-1})$ for $\vphi$ in a countable dense subset of the nuclear Fr\'echet space $(\mS,
\N{\cdot}_{p}, p\geq 1)$.
\end{rem}
Thanks to Mitoma's result, it suffices to have a tightness criterion in $\R{}$. We recall here the usual
result (cf. Billingsley \cite{cf:BillConvProb}):
A sequence of $(\Om^{N}, \mF_{t}^{N})$-adapted processes $(Y^{N})$ with paths in $\mC([0,T], \R{})$ is tight
if both of the following conditions hold:
\begin{itemize}
	\item Condition [T]: for all $t\leq T$ and $\delta>0$, there exists $C>0$ such that
	\begin{equation}
\label{eq:condT}
\tag{$T_{t, \delta, C}$}
\sup_{N} \bP\left(|Y^{N}_{t}|>C \right)\leq \delta,
\end{equation}
	\item Condition [A]: for all $\eta_{1},\eta_{2}>0$, there exists $C>0$ and $N_{0}$ such that for all
$\mF^{N}$-stopping times $\tau_{N}$,\begin{equation}
\label{eq:condA}
\tag{$A_{\eta_{1}, \eta_{2}, C}$}
\sup_{N\geq N_{0}}\sup_{\theta\leq C}\bP\left(\left|Y^{N}_{\tau_{N}} - Y^{N}_{\tau_{N}+\theta}\right|\geq
\eta_{2}\right)\leq \eta_{1}.
\end{equation}
\end{itemize}

\begin{proof}[Proof of Theorem \ref{theo:tightthetaxi}]
\begin{enumerate}
	\item Tightness of $(M^{N, (\om)})$: for a fixed realization of the disorder $(\om)$, for fixed
$\vphi\in\mS$, we have:
\begin{itemize}
	\item For all $t\in[0, T]$, for all $\delta>0$, for all $C>0$, 
\begin{align*}
	\bP\left(\left|M^{N, (\om)}_{t}(\vphi)\right|>C\right)&\leq \frac{\bE\left[\sup_{t\leq
T}\left\{M_{t}^{N, (\om)}(\vphi)^{2}\right\}\right]}{C^{2}},\\
	&\leq \frac{\bE\left[\sup_{t\leq T}\N{M_{t}^{N, (\om)}}_{-3,2\alpha}^{2}\N{\vphi}_{3,
2\alpha}^{2}\right]}{C^{2}},\\
	&\leq \frac{C\N{\vphi}_{3, 2\alpha}^{2}}{a^{2}}\sup_{N}\frac{1}{N} \Som{i}{1}{N}{\left(1+
|\om_{i}|^{4\alpha}\right)}, \text{(cf. \eqref{eq:remsupNMt})},\\
	&\leq \delta,
\end{align*}
for a suitable $C>0$ (depending on $(\om)$). Condition [T] is proved.
\item Let us verify Condition [A]: 
For every $\vphi\in\mS$, for every $\delta, \theta, \eta_{1}, \eta_{2}>0, \theta\leq \delta$, for every
stopping time $\tau_{N}$,
\begin{align*}
u_{N}	&:=\bP\left(\left| M^{N}_{\tau_{N}+\theta}(\vphi)- M^{N}_{\tau_{N}}(\vphi)\right|>\eta_{2}\right)\leq
\frac{1}{\eta_{2}^{2}} \bE\left[\left| M^{N}_{\tau_{N}+\theta}(\vphi)-
M^{N}_{\tau_{N}}(\vphi)\right|^{2}\right],\\
&\leq\frac{1}{\eta_{2}^{2}}\bE\left[\int_{\tau_{N}}^{\tau_{N}+\theta}{\cro{\nu_{s}^{N}}{\vphi'(y, \pi)^{2}}\d
s}\right],\\
&\leq
\N{\vphi}_{6,\alpha}^{2}\frac{1}{\eta_{2}^{2}}\bE\left[\int_{\tau_{N}}^{\tau_{N}+\theta}{\int_{\Sun\times\R{}}
{\N{H_{y,\pi}}_{-6, \alpha}^{2}\d\nu_{s}^{N}}\d s}\right],\\
&\leq
\N{\vphi}_{6,\alpha}^{2}\frac{C}{\eta_{2}^{2}}\bE\left[\int_{\tau_{N}}^{\tau_{N}+\theta}{\frac{1}{N}\Som{i}{1}
{N}{(1+|\om_{i}|^{4\alpha})}\d s}\right],\ \text{(cf. \eqref{eq:injsobo3} and \eqref{eq:NHxo}),}\\
	&\leq
\N{\vphi}_{6,\alpha}^{2}\frac{C\delta}{\eta_{2}^{2}}\sup_{N}\left(\frac{1}{N}\Som{i}{1}{N}{(1+|\om_{i}|^{
4\alpha})}\right).
\end{align*}
This last term is lower or equal than $\eta_{1}$ for $\delta$ sufficiently small (depending on $(\om)$).
\end{itemize}	

\item Tightness of $(\Theta_{N})$: we need to be more careful here, since the tightness is \emph{in law w.r.t.
the disorder}. 
Let $(\vphi_{j})_{j\geq 1}$ be a countable family in the nuclear Fr\'echet space $\mS$. Without any
restriction, we can always suppose that $\N{\phi_{j}}_{6, \alpha}=1$, for every $j\geq 1$. We define the
following decreasing sequences (indexed by $J\geq 1$) of subsets of $\mM_{1}\left(\mC([0, T], \mS')\right)$:
\begin{align*}
	K_{1}^{\eps}(\vphi_{1}, \dots, \vphi_{J})&:=\ens{P}{\forall t, \forall 1\leq j\leq J,\
P\Pi_{\vphi_{j}}^{-1}\text{ satisfies $(T_{t, \delta, C_{1}})$}},\\
	K_{2}^{\eps}(\vphi_{1}, \dots, \vphi_{J})&:= \ens{P}{\forall 1\leq j\leq J, \forall
\eta_{1},\eta_{2}>0, P\Pi_{\vphi_{j}}^{-1} \text{ satisfies $(A_{\eta_{1}, \eta_{2}, C_{2}})$}},
\end{align*}
where $C_{1}=C_{1}(\eps, \delta)$, $C_{2}=C_{2}(\eps, \eta_{1}, \eta_{2})$ will be precised later.
By construction and by Mitoma's theorem (cf. Remark \ref{rem:phisep}), \[K^{\eps}:=\bigcap_{J}
(K_{1}^{\eps}(\vphi_{1}, \dots, \vphi_{J})\cap K_{2}^{\eps}(\vphi_{1}, \dots, \vphi_{J}))\] is a relatively
compact subset of $\mM_{1}\left(\mC([0, T], \mS')\right)$. In order to prove tightness of $(\Theta_{N})$, it
is sufficient to prove that, for all $\eps>0$, \begin{equation}
\label{eq:tightnessTheta}
\forall i=1,2,\ \limsup_{N}\Theta_{N}\left(\bigcup_{J}K_{i}^{\eps}(\phi_{1}, \dots, \phi_{J})^{c}\right)\leq
\eps.
\end{equation}
For $\eps>0$, let $A=A(\eps)$ such that $\liminf_{N\tend\infty}\bbP\left(A_{N}\leq A\right)\geq 1-\eps,$ and
\[\liminf_{N\tend\infty}\bbP\left(\frac{1}{N}\Som{i}{1}{N}{(1+|\om_{i}|^{4\alpha})} + A_{N}(\om_{1},\dots,
\om_{N}) \leq A\right)\geq 1-\eps,\] where $A_{N}$ is the random variable defined in Proposition
\ref{prop:estimeta1}. We define the corresponding constants (for a sufficiently large constant $C$):
\[C_{1}(\eps, \delta):= \sqrt{\frac{A(\eps)}{\delta}},\quad C_{2}(\eps, \eta_{1}, \eta_{2}):=
\frac{\eta_{1}\eta_{2}^{2}}{CA(\eps)}.\]Then,
\begin{align*}
	\Theta_{N}(K_{1}^{\eps}(\phi_{1}, \dots, \phi_{J})) &\geq \bbP\left((\om),\ \forall t,\ \forall 1\leq
j\leq J,\ \forall \delta,\ \frac{\bE\left[\left|\cro{\eta_{t}^{N,
(\om)}}{\phi_{j}}\right|^{2}\right]}{C_{1}(\delta, \eps)^{2}}\leq \delta\right),\\
	&\geq \bbP\left((\om),\ \sup_{t\leq T}\bE\left[\N{\eta_{t}^{N, (\om)}}_{-6,\alpha}^{2}\right]\leq
A\right),\ \text{(by definition of $C_{1}$)},\\
	&\geq \bbP\left(A_{N}\leq A(\eps)\right),\ \text{(cf. \eqref{eq:injsobo3} and \eqref{eq:estimeta1})}.
\end{align*}Letting $J\tend\infty$ in the latter inequality, we obtain:
$\Theta_{N}(\bigcup_{J}K_{1}^{\eps}(\phi_{1}, \dots, \phi_{J})^{c})\leq \bbP(A_{N}>A).$ Taking on both sides
$\limsup_{N\tend\infty}$, we get the result.

Furthermore, for $\eta_{2}>0$, $0<\theta\leq C_{2}$ and $\tau_{N}\leq T$ a stopping time, for all $1\leq j\leq
J$,
\begin{align*}
	\bP\left(\left|\int_{\tau_{N}}^{\tau_{N}+\theta}{\cro{\eta_{s}^{N,
(\om)}}{\mL_{s}^{\nu^{N}}(\vphi_{j})}\d s}\right|\geq
\eta_{2}\right)&\leq\frac{1}{\eta_{2}^{2}}\bE\left[\left|\int_{\tau_{N}}^{\tau_{N}+\theta}{\cro{\eta_{s}^{N,
(\om)}}{\mL_{s}^{\nu^{N}}(\vphi_{j})}\d s}\right|^{2}\right],\\
	&\leq \frac{C_{2}}{\eta_{2}^{2}} \bE\left[\int_{\tau_{N}}^{\tau_{N}+\theta}{\left|\cro{\eta_{s}^{N,
(\om)}}{\mL_{s}^{\nu^{N}}(\vphi_{j})}\right|^{2}\d s}\right],\\
	&\leq \frac{C_{2}}{\eta_{2}^{2}} \int_{0}^{T}{\bE\left|\cro{\eta_{s}^{N,
(\om)}}{\mL_{s}^{\nu^{N}}(\vphi_{j})}\right|^{2}\d s},\\
	&\leq \frac{CC_{2}}{\eta_{2}^{2}}\int_{0}^{T}{\bE\left[\N{\eta_{s}^{N}}_{-3, 2\alpha}^{2}\right]\d
s},\\
	&\leq \frac{CTC_{2}}{\eta_{2}^{2}}A_{N},\ \text{(cf. \eqref{eq:estimeta1})}.
\end{align*}And,\[\bP\left(\left| M^{N}_{\tau_{N}+\theta}(\vphi_{j})-
M^{N}_{\tau_{N}}(\vphi_{j})\right|>\eta_{2}\right)\leq
\frac{CC_{2}}{\eta_{2}^{2}}\left(\frac{1}{N}\Som{i}{1}{N}{(1+|\om_{i}|^{4\alpha})}\right).\]
So, for all $j\geq 1$, by definition of $C_{2}$, \[\bP\left(\left|\eta_{\tau_{N}+\theta}^{N, (\om)}(\vphi_{j})
- \eta_{\tau_{N}}^{N, (\om)}(\vphi_{j})\right|\geq \eta_{2}\right)\leq \frac{\eta_{1}}{A(\eps)}\left(A_{N} +
\frac{1}{N}\Som{i}{1}{N}{(1+|\om_{i}|^{4\alpha})}\right).\]
Consequently, \[\Theta_{N}(K_{2}^{\eps}(\phi_{1}, \dots, \vphi_{J})) \geq \bbP\left(A_{N} +
\frac{1}{N}\Som{i}{1}{N}{(1+|\om_{i}|^{4\alpha})} > A(\eps)\right).\]Letting $J\tend\infty$, we get
$\limsup_{N}\Theta_{N}\left(\bigcup_{J}K_{2}^{\eps}(\vphi_{1}, \dots, \vphi_{J})^{c}\right)\leq \eps$. Eq.
\eqref{eq:tightnessTheta} is proved.\qedhere
\end{enumerate}
\end{proof}
\subsubsection{Identification of the limit}
The proof of the fluctuations result will be complete when we identify any possible limit.
\begin{prop}[(Identification of the initial value)]
\label{prop:idiniteta}
The random variable $(\om)\mapsto \mL\left(\eta_{0}^{N, (\om)}\right)$ converges in law to the random variable
$\om\mapsto \mL(X(\om))$, where for all $\om$, $X(\om)=C(\om) +Y$, with $Y$ a centered Gaussian process with
covariance $\Gamma_{1}$. Moreover $\om\mapsto C(\om)$ is a Gaussian process with covariance $\Gamma_{2}$,
where $\Gamma_1$ and $\Gamma_2$ are defined in \eqref{eq:covC1} and \eqref{eq:covC2}.
\end{prop}

\begin{proof}
For simplicity, we only identify here the law of $\cro{\eta_{0}^{N, (\om)}}{\vphi}$ for all $\vphi$. The same
proof works for the law of finite-dimensional distributions $\left(\cro{\eta_{0}^{N, (\om)}}{\vphi_{1}},
\dots,\cro{\eta_{0}^{N, (\om)}}{\vphi_{p}}\right)$, $p\geq 1$. We write $\Gamma_{i}$ for $\Gamma_{i}(\vphi,
\vphi)$, $i=1,2$. One has:
\begin{align*}
\cro{\eta_{0}^{N, (\om)}}{\vphi} &= \frac{1}{\sqrt{N}}\Som{i}{1}{N}{\left(\vphi(\xi^{i}, \om_{i}) -
\int_{\Sun}{\vphi(x, \om_{i})\lambda(\d x)}\right)}\\ &+
\frac{1}{\sqrt{N}}\Som{i}{1}{N}{\left(\int_{\Sun}{\vphi(x, \om_{i})\lambda(\d x)} -
\cro{\nu_{0}}{\vphi}\right)},\\
&=:A^{N, (\om)} + B^{N, (\om)}.
\end{align*}
It is easy to see that $B^{N, (\om)}$ converges in law to $Z_{2}\sim \mN(0, \Gamma_{2})$. Moreover, for
$\bbP$-almost
every $(\om)$, $A^{N, (\om)}$ converges in law to $Z_{1}\sim\mN(0, \Gamma_{1})$ (see Billingsley
\cite{cf:Billmeasure}, Th.~27.3 p. 362). That means that for all $u\in\R{}$,
$\psi_{A^N}(u):=\bE_{\lambda}\left(e^{iu A^{N, (\om)}}\right)$ converges to
$\psi_{Y}(u):=e^{-\frac{u^2}{2\Gamma_{1}}}$.
But, then, for all $F\in\mC_{b}(\R{})$,
\begin{align*}
 \bbE\left[F\left(\bE_{\lambda}\left[e^{iu \cro{\eta_{0}^{N, (\om)}}{\vphi}}\right]\right)\right] &=
\bbE\left[F\left(\bE_{\lambda}\left[e^{iu (A^{N, (\om)}+B^{N,
(\om)})}\right]\right)\right]=\bbE\left[F\left(e^{iu
B^{N, (\om)}}\psi_{N}(u)\right)\right].
\end{align*}
Since $\psi_{N}(u)$ converges almost surely to a constant, the limit of the expression above exists (Slutsky's
theorem) and is equal to $\bbE\left[F\left(e^{iuZ_{2} - \frac{u^2}{2\Gamma_{1}}}\right)\right]$.  
\end{proof}
\begin{prop}[(Identification of the martingale part)]
\label{prop:idWeta}
For $\bbP$-almost every $(\om)$, the sequence $(M^{N, (\om)})$ converges in law in $\mC([0, T], \mS')$ to a
Gaussian process $W$ with covariance defined in \eqref{eq:Wmart}.
\end{prop}
\begin{proof}
For fixed $(\om)$, $(M^{N, (\om)})$ is a sequence of uniformly square-integrable continuous martingales (cf.
Remark \ref{rem:remsupNMt}), which is tight in $\mC([0, T], \mS')$. Let $W_{1}$ and $W_{2}$ be two
accumulation points (continuous square-integrable martingales which \emph{a priori} depend on $(\om)$) and
$(M^{\phi(N), (\om)})$ and $(M^{\psi(N), (\om)})$ be two subsequences converging to $W_{1}$ and $W_{2}$,
respectively. Note that we can suppose that $\phi(N)\leq \psi(N)$ for all $N$. For all $\vphi\in \mS$,
$\lim_{N\tend\infty} \cro{M^{\phi(N), (\om)}(\vphi)}{M^{\psi(N), (\om)}(\vphi)}_{t} =
\cro{W_{1}(\vphi)}{W_{2}(\vphi)}_{t}$, for all $t$, and \[\cro{M^{\phi(N), (\om)}(\vphi)}{M^{\psi(N),
(\om)}(\vphi)}_{t} = \int_{0}^{t}{\cro{\nu_{s}^{\phi(N)}}{(\vphi')^{2}}\d s}.\] We now have to identify the
limit: we already know that for $\bbP$-almost every realization of the disorder $(\om)$, $(\nu^{N , (\om)})$
converges in law to $P$. But, the latter expression, seen as a function of $\nu$, is continuous. So
$\cro{W_{1}(\vphi)}{W_{2}(\vphi)}_{t}=\int_{0}^{t}{\cro{P_{s}}{(\vphi')^{2}}}.$ So $W_{1}-W_{2}$ is a
continuous square integrable martingale whose Doob-Meyer process is $0$. So $W_{1}=W_{2}$ and is characterized
as the Gaussian process with covariance given in \eqref{eq:Wmart}. The convergence follows.
\end{proof}

\begin{proof}[Proof of the independence of $W$ and $X$]
We prove more : the triple $(Y, C, W)$ is independent.
For sake of simplicity, we only consider the case 
of $(Y(\vphi), C(\vphi), W_{t}(\vphi))$ for
fixed $t$ and $\vphi$. 

Let us first recall some notations: let $A^{N, (\om)}$, $B^{N, (\om)}$ and $M^{N,
(\om)}_{t}(\vphi)$ be the random variables defined in the proof of Proposition \ref{prop:idiniteta} and
\ref{prop:idWeta} and let $\psi_{A^N}(u):=\bE\left( e^{iu A^{N, (\om)}} \right)$, $\psi_{B^N}(v):= \bE\left(
e^{iv B^{N, (\om)}} \right)$, $\psi_{M^N}(w):=\bE\left( e^{iw M^{N, (\om)}_{t}(\vphi)} \right)$ be their
characteristic functions ($u, v, w \in \R{}$).
We know that, for almost every $(\om)$, $\psi_{A^N}(u)$ converges to $\psi_Y(u)= e^{-\frac{u^2}{2\Gamma_{1}}}$
and that $\psi_{M^N}(w)$ converges to the deterministic function $\psi_W(w):=
\bE\left(e^{iw
W_{t}(\vphi)}\right)$. But, if $\psi_C(v)=\bbE\left( e^{iw C} \right)$, then, for all $u, v, w\in\R{}$, using
the
independence of the Brownian with the initial conditions,
\[\bbE\left(\bE\left(e^{iuA^{N,
(\om)} + iv B^{N, (\om)} +iw M^{N, (\om)}_{t}(\vphi)} \right)
- e^{iv B^{N, (\om)}} \psi_{A^N}(u)\psi_{M^N}(w)\right)=0.\] Using Slutsky's theorem, we see that any limit
couple $(Y, C, W)$ satisfies \[\bbE\left(
\bE\left(e^{iu Y+ ivC +iw W_{t}(\vphi)}\right)
\right)=\psi_{Y}(u) \psi_{C}(v) \psi_{W}(w).\]  which is the
desired result.
\end{proof}

We recall that the limit second order differential operator $\mL_{s}$ is defined by
\[\mL_{s}(\vphi)(y, \pi):=  \frac{1}{2} \vphi''(y, \pi) + \vphi'(y, \pi)(b[y, P_{s}] + c(y, \pi)) +
\cro{P_{s}}{\vphi'(\cdot,\cdot)b(\cdot, y, \pi)}.\]
As in Lemma \ref{lem:estimLmu}, we can prove the following:
\begin{lem}
\label{lem:estimmL}
Assume \eqref{eq:Hbc}. Then for every $N$, $s\leq T$, $(\om)$, the operator $\mL_{s}$ and $\mL_{s}^{\nu^{N}}$
are linear continuous from $\mS$ to $\mS$ and
\[\N{\mL_{s}(\vphi)}_{6, \alpha} \leq C \N{\vphi}_{8, \alpha},\]
\[\N{\mL_{s}^{\nu^{N}}(\vphi)}_{6, \alpha} \leq C \N{\vphi}_{8, \alpha}.\]
\end{lem}

We are now in position to prove Theorem \ref{theo:fluctquenchedinit}:
\begin{proof}[Proof of Theorem \ref{theo:fluctquenchedinit}]
Let $\Theta$ be an accumulation point of $\Theta_{N}$. Thus, for a certain subsequence (which will be also
denoted as $N$ for notations purpose), the random variable $(\om)\mapsto \mH^{N, (\om)}$ converges in law to a
random variable $\mH$ with values in $\mM_{1}(\mC([0, T], \mS'))$ with law $\Theta$. Applying Skorohod's
representation theorem, there exists some probability space $(\Om^{(1)}, \bP^{(1)}, \mF^{(1)})$ and random
variables defined on $\Om^{(1)}$, $\om_{1}\mapsto H^{N}(\om_{1})$ and $\om_{1}\mapsto H(\om_{1})$ such that
$H^{N}$ has the same law as $(\om)\mapsto \mH^{N, (\om)}$, $H$ has the same law as $\mH$, and for
$\bP^{(1)}$-almost every $\om_{1}\in \Om^{(1)}$, $H^{N}(\om_{1})$ converges to $H(\om_{1})$ in
$\mM_{1}(\mC([0, T], \mS'))$. 

An easy application of Proposition \ref{prop:estimsupetat} and Borel-Cantelli's Lemma shows that
$\bP^{(1)}$-almost surely, $\bE\left(\sup_{t\leq T} \N{\eta_{t}^{\om_{1}}}_{-6, \alpha}\right)<\infty.$ Then
we know from Lemma \ref{lem:estimmL} that the integral term $\int_{0}^{t}{\mL_{s}^{*}\eta_{s}^{\om_{1}}\d s}$
makes sense as a Bochner's integral in $W_{0}^{-8, \alpha}\subseteq \mS'$.

Let $\eta^{N, \om_{1}}$ with law $H^{N}(\om_{1})$; $\eta^{N, \om_{1}}$ converges in law to some
$\eta^{\om_{1}}$ with law $H(\om_{1})$. By uniqueness in law convergence, using Propositions
\ref{prop:idiniteta} and \ref{prop:idWeta}, we see that $(\eta^{\om_{1}}_{0}, W)$ as the same law as
$(X(\om_{1}), W)$. For fixed $\vphi \in
\mS$, we define $F_{\vphi}$ from $\mC([0, T], \mS')$ into $\R{}$ by $F_{\vphi}(\gamma):=
\cro{\gamma_{t}}{\vphi} - \cro{\gamma_{0}}{\vphi} - \int_{0}^{t}{\cro{\gamma_{s}}{\mL_{s}\vphi}\d s}$.
The function $F_{\vphi}$ is continuous and since $\eta^{N, \om_{1}}$ converges in law to $\eta^{\om_{1}}$, the
sequence $(F_{\vphi}(\eta^{N, \om_{1}}))$ converges in law to $F_{\vphi}(\eta^{\om_{1}})$.
To prove the theorem, it remains to show that the law of the term $\int_{0}^{t}{\cro{\eta^{N,
\om_{1}}_{s}}{\mL_{s}^{\nu^{N}}\vphi - \mL_{s}\vphi}\d s}$ converges in law to $0$. We show that there is
convergence in probability:
For all $\eps>0$, for all $A>0$, using Proposition \ref{prop:estimsupetat}, Lemma \ref{lem:estimmL}, and
Cauchy-Schwarz's inequality, 
\begin{align*}
	U_{N, \eps} &:= \bP^{(1)}\left(\bE\left[\int_{0}^{t}{\left|\cro{\eta^{N,
\om_{1}}_{s}}{(\mL_{s}^{\nu^{N}} - \mL_{s})(\vphi)}\right|\d s}\right]>\eps\right),\\
	&= \bbP\left(\bE\left[\int_{0}^{t}{\left|\cro{\eta^{N, (\om)}_{s}}{(\mL_{s}^{\nu^{N}} -
\mL_{s})(\vphi)}\right|\d s}\right]>\eps\right),\\
	&\leq  \bbP\left(\int_{0}^{t}{\bE\left[\N{\eta^{N, (\om)}_{s}}_{-6,
\alpha}^{2}\right]^{1/2}\bE\left[\N{(\mL_{s}^{\nu^{N}} - \mL_{s})(\vphi)}_{6, \alpha}^{2}\right]^{1/2}\d
s}>\eps\right),\\
	&\leq \bbP\left(C_{N}(\om_{1}, \dots, \om_{N})^{1/2}\int_{0}^{t}{\bE\left[\N{(\mL_{s}^{\nu^{N}} -
\mL_{s})(\vphi)}_{6, \alpha}^{2}\right]^{1/2}\d s}>\eps\right)\ \text{(cf. Prop \ref{prop:estimsupetat})},\\
	&\leq \bbP\left(\int_{0}^{t}{\bE\left[\N{(\mL_{s}^{\nu^{N}} - \mL_{s})(\vphi)}_{6,
\alpha}^{2}\right]^{1/2}\d s}>\frac{\eps}{\sqrt{A}}\right) + \bbP\left(C_{N} > A\right).
\end{align*}
Using \eqref{eq:estimsupetatCN}, it suffices to prove that, for all $\eps>0$,
\begin{equation}
\label{eq:estimfinalconvloi}
\limsup_{N\tend\infty} \bbP\left(\int_{0}^{t}{\bE\left[\N{(\mL_{s}^{\nu^{N}} - \mL_{s})(\vphi)}_{6,
\alpha}^{2}\right]^{1/2}\d s}>\eps\right)= 0.
\end{equation} 
Indeed, for every $\vphi\in \mS$,
\begin{align*}
	\mU^{N}_{s}(\vphi)(y, \pi) &:= (\mL_{s}^{\nu^{N}} - \mL_{s})(\vphi)(y, \pi)= \vphi'(y, \pi) (b[y,
\nu_{s}^{N}] - b[y, P_{s}]).
\end{align*}
An analogous calculation as in Lemma \ref{lem:estimLmu} shows that, using Lipschitz assumptions on  $b$, and
Proposition \ref{prop:estimpnonlin}:
\[\bE\left[\N{\sup_{s\leq t}\mU^{N}_{s}(\vphi)}_{6, \alpha}^{2}\right] \leq \N{\vphi}_{8, \alpha}^{2}(C/N +
D_{N}(\om_{1}, \dots, \om_{N})),\]
with the property that $\lim_{A\tend\infty}\limsup_{N}\bbP(N D_{N}>A)= 0$. Equation
\eqref{eq:estimfinalconvloi} is a direct consequence.

Since there is uniqueness in law in \eqref{eq:OUeta}, $\Theta$ is perfectly defined, and thus, unique. The
convergence follows.
\end{proof}

\section{Proofs for the fluctuations of the order parameters}
\label{sec:prooffluctrpsi}
We end by the proofs of paragraph \ref{sec:fluctrpsi}.

\subsection{Proof of Proposition \ref{prop:convfluctr}}
\begin{enumerate}
	\item This is straightforward since $r^{N, (\om)} = \left|\cro{\nu^{N, (\om)}}{e^{ix}}\right|$ and
since for $\bbP$-almost every disorder $(\om)$, $\nu^{N, (\om)}$ converges weakly to $P$.
	\item The following sequences are well defined:
$\forall k\geq 0,$
\begin{align*}
u_{k}(t)&:= \int_{\Sun\times\R{}}{e^{-|\om|}\om^{k}\cos(\theta)P_{t}(\d\theta, \d\om)},\\
v_{k}(t)&:= \int_{\Sun\times\R{}}{e^{-|\om|}\om^{k}\sin(\theta)P_{t}(\d\theta, \d\om)}.
\end{align*}
Let $E=(\ell_{\infty}(\bN), \N{\cdot}_{\infty})$ be the Banach space of real bounded sequences endowed with
its usual $\N{\cdot}_{\infty}$ norm, ($\N{u}_{\infty}= \sup_{k\geq 0} |u_{k}|$). For all $t>0$, let $\mA_{t}:
E\times E \tend E\times E$, be the following linear operator (where $(u, v)$ is a typical element of $E\times
E$):
For all $k\geq 0$
\[\left\{\begin{array}{ccc}
	\mA_{t}(u,0)_{k} &=& -\frac{1}{2} u_{k} - \alpha_{k}(t) v_{0} + \beta_{k}(t) u_{0} -K v_{k+1},\\
	\mA_{t}(0,v)_{k} &=& -\frac{1}{2} v_{k} + \gamma_{k}(t) v_{0} - \alpha_{k}(t) u_{0} +K u_{k+1},
\end{array}\right.\]
where,\begin{align*}
\alpha_{k}(t)&=\cro{P_{t}}{e^{-|\om|}\om^{k}\cos(\cdot)\sin(\cdot)},\\
\beta_{k}(t)&=\cro{P_{t}}{e^{-|\om|}\om^{k}\sin^{2}(\cdot)},\\
\gamma_{k}(t)&=\cro{P_{t}}{e^{-|\om|}\om^{k}\cos^{2}(\cdot)}.
\end{align*}
$(t,u,v)\mapsto \mA_{t}\cdot (u,v)$ is globally Lipschitz-continuous map from $[0,T]\times E\times E$ into
$E\times E$ and one easily verifies considering \eqref{eq:weakMKV} (in the case of the sine-model) and
developing the sine interaction that $t\mapsto(u(t), v(t))$ satisfies in $E\times E$ the following linear
non-homogeneous Cauchy Problem:
\[
\left\{\begin{array}{ccl}
\frac{\d}{\d t} (u(t), v(t)) &=& \mA_{t}\cdot(u(t), v(t)),\\
 u_{k}(0) &=&\cro{P_0}{e^{-|\om|}\om^{k}\cos(\cdot)},\ \forall k\geq0\\
 v_{k}(0)&=& \cro{P_0}{e^{-|\om|}\om^{k}\sin(\cdot)},\ \forall k\geq0.	
\end{array}\right.\]
Let us suppose that there exists some $t_{0}\in[0,T]$ such that $r_{t_{0}}=0$, namely $u_{0}(t_{0})=
v_{0}(t_{0}) =0$. Then, if $(\tilde{u}, \tilde{v})$ is the constant function on $[0, T]$ such that for all
$k\geq 0$, $\tilde{u}_{k} \equiv u_{k}(t_{0}), \tilde{v}_{k} \equiv v_{k}(t_{0})$, then $(\tilde{u},
\tilde{v})$ satisfy the same Cauchy Problem as $(u, v)$ with initial condition at time $t_{0}$. By
Cauchy-Lipschitz theorem, both functions coincide on $[0, T]$. In particular, $u_{0}$ and $v_{0}$ are always
zero and thus $r\equiv0$.
	\item We suppose \eqref{eq:rtpos}.
	A simple calculation shows that the fluctuations process $\mR^{N}$ verifies for all $t\in[0, T]$,
	\begin{align*}
	\mR_{t}^{N,(\om)} &= \frac{\cro{\eta_{t}^{N,(\om)}}{\cos(\cdot)} \cro{\nu_{t}^{N, (\om)}+
P_t}{\cos(\cdot)} + \cro{\eta_{t}^{N,(\om)}}{\sin(\cdot)} \cro{\nu_{t}^{N, (\om)}+
P_t}{\sin(\cdot)}}{r_{t}^{N, (\om)}+ r_{t}},\\
	&= \frac{\Re\left(\cro{\eta_{t}^{N, (\om)}}{e^{ix}}\overline{\cro{\nu_{t}^{N,
(\om)}+P_t}{e^{ix}}}\right)}{\left|\cro{\nu_{t}^{N, (\om)}}{e^{ix}}\right|+ r_t}.
	\end{align*}
	Let $u^{N, (\om)}:= \cro{\nu^{N, (\om)}}{e^{ix}}$, $v^{N, (\om)}:= \cro{\eta^{N, (\om)}}{e^{ix}}$ and
$u:= \cro{P}{e^{ix}}$, $v^{\om}:= \cro{\eta^{\om}}{e^{ix}}$ be their corresponding limits. The result follows
if we prove the following property: the random variables $(\om)\mapsto \mL\left(u^{N, (\om)}, v^{N,
(\om)}\right)$ converges in law to the random variable $\om\mapsto \mL\left(u , v^{\om}\right)$. The tightness
of this random variable follows from the convergence of both empirical measure and fluctuations process. As
already said in Remark \ref{rem:Rtp}, it suffices to prove the convergence of the finite-dimensional marginals
$(\sou{u}^{N, (\om)}_{\sou{t}}, \sou{v}^{N, (\om)}_{\sou{t}})=\left((u^{N, (\om)}_{t_{1}},\dots, u^{N,
(\om)}_{t_{p}}), (v^{N, (\om)}_{t_{1}}, \dots, v^{N, (\om)}_{t_{p}})\right)$, for all element of $[0, T]$,
$t_{1}, \dots, t_{p}$, $p\geq 1$.
	
Since the limit of $(u^{N, (\om)})$ is a constant, this is mainly a consequence of Slutsky's theorem. But
since this is a convergence \emph{in law with respect to the disorder}, one has to adapt the proof. We prove
the following: $\forall G\in\mC_{b}^{1}(\R{})$, $\forall \sou{r}=\left(r_{1}, \dots, r_{p}\right)\in\R{p}$,
$\forall \sou{s}=\left(s_{1}, \dots, s_{p}\right)\in\R{p}$,
\begin{align*}
\bbE\left[G\left(\vphi_{(\sou{u}^{N, (\om)}_{\sou{t}}, \sou{v}^{N,
(\om)}_{\sou{t}})}(\sou{r},\sou{s})\right)\right]
\underset{N\tend\infty}{\tend} \bbE\left[G\left(\vphi_{(\sou{u}_{\sou{t}},
\sou{v}_{\sou{t}}^{\om})}(\sou{r},\sou{s})\right)\right],
\label{eq:Hphi}
\end{align*}
where $\vphi_{(\sou{X}, \sou{Y})}(\sou{r},\sou{s}) = \bE\left[e^{i \sou{r}\cdot \sou{X}+i \sou{s}\cdot
\sou{Y}}\right]$ is the characteristic function of the couple $(\sou{X}, \sou{Y})$. Indeed, we have
successively:
\begin{align*}
a_N&:=\left|\bbE\left[G\left(\vphi_{(\sou{u}^{N, (\om)}_{\sou{t}}, \sou{v}^{N,
(\om)}_{\sou{t}})}(\sou{r},\sou{s})\right)\right] - \bbE\left[G\left(\vphi_{(\sou{u}_{\sou{t}},
\sou{v}^{\om}_{\sou{t}})}(\sou{r},\sou{s})\right)\right]\right|,\\
&\leq \left|\bbE\left[G\left(\vphi_{(\sou{u}^{N, (\om)}_{\sou{t}}, \sou{v}^{N,
(\om)}_{\sou{t}})}(\sou{r},\sou{s})\right)\right] - \bbE\left[G\left(\vphi_{(\sou{u}_{\sou{t}}, \sou{v}^{N,
(\om)}_{\sou{t}})}(\sou{r}, \sou{s})\right)\right]\right|\\
&+ \left|\bbE\left[G\left(\vphi_{(\sou{u}_{\sou{t}}, \sou{v}^{N,
(\om)}_{\sou{t}})}(\sou{r},\sou{s})\right)\right] - \bbE\left[G\left(\vphi_{(\sou{u}_{\sou{t}},
\sou{v}^{\om}_{\sou{t}})}(\sou{r},\sou{s})\right)\right]\right|,\\
&\leq C\bbE\left|\vphi_{(\sou{u}^{N, (\om)}_{\sou{t}}, \sou{v}^{N, (\om)}_{\sou{t}})}(\sou{r},\sou{s}) -
\vphi_{(\sou{u}_{\sou{t}}, \sou{v}^{N, (\om)}_{\sou{t}})}(\sou{r},\sou{s})\right|\\
&+ \left|\bbE\left[G\left(\vphi_{(\sou{u}_{\sou{t}}, \sou{v}^{N,
(\om)}_{\sou{t}})}(\sou{r},\sou{s})\right)\right] - \bbE\left[G\left(\vphi_{(\sou{u}_{\sou{t}},
\sou{v}^{\om}_{\sou{t}})}(\sou{r},\sou{s})\right)\right]\right|,\\
&\leq  C\bbE\bE\left|e^{i\sou{r}\cdot\sou{u}^{N, (\om)}_{\sou{t}}}-e^{i\sou{r}\cdot\sou{u}_{\sou{t}}}\right|\\
&+ \left|\bbE\left[G\left(\vphi_{(\sou{u}_{\sou{t}}, \sou{v}^{N,
(\om)}_{\sou{t}})}(\sou{r},\sou{s})\right)\right] - \bbE\left[G\left(\vphi_{(\sou{u}_{\sou{t}},
\sou{v}^{\om}_{\sou{t}})}(\sou{r},\sou{s})\right)\right]\right|.\\
\end{align*}
But, we have $\bE\left|e^{i\sou{r}_{\sou{t}}\cdot\sou{u}^{N,
(\om)}_{\sou{t}}}-e^{i\sou{r}\cdot\sou{u}_{\sou{t}}}\right|\leq \min\left(2, |\sou{r}||\sou{u}^{N,
(\om)}_{\sou{t}}-\sou{u}_{\sou{t}}|\right)$. So, for all $\eps>0$, \[\bE\left|e^{i\sou{r}\cdot\sou{u}^{N,
(\om)}_{\sou{t}}}-e^{i\sou{r}\cdot\sou{u}_{\sou{t}}}\right|\leq \eps|\sou{r}| + 2\bP\left(\left|\sou{u}^{N,
(\om)}_{\sou{t}}-\sou{u}_{\sou{t}}\right|>\eps\right).\] Taking $\limsup_{N\tend\infty}$, and letting
$\eps\tend 0$, we get $\lim a_{N}=0$. The result follows.\qedhere
\end{enumerate}

\subsection{Proof of Proposition \ref{prop:convfluctpsi}}
The proof is similar to the previous one and relies on the two following equalities :
\begin{align*}
 \zeta^{N, (\om)} &= \frac{\cro{P}{e^{ix}}}{r^{N, (\om)}},\\
 \sqrt{N}\left(\zeta^{N, (\om)}- \zeta\right) &= \frac{1}{r\cdot r^{N, (\om)}} \left(r\cro{\eta^{N,
(\om)}}{e^{ix}} + \cro{P}{e^{ix}}\mR^{N, (\om)}\right).
\end{align*}
\appendix

\section{Proof of Proposition \ref{prop:estimpnonlin}}
\label{app:propZN}
Thanks to the Lipschitz continuity of $b$ and $c$, introducing $\nu$ as the empirical measure corresponding to
$(x^{\om_{i}}, \om_{i})$, we have, (inserting $b[x_{s}^{\om_{i}}, \nu_{s}^{N}] - b[x_{s}^{\om_{i}}, \nu_{s}]$
in the $b$ term), 
\begin{align*}
	\bE\left[\sup_{s\leq t}\left|x_{s}^{i, N} - x_{s}^{\om_{i}}\right|^{2}\right] &\leq C
\left(\int_{0}^{t}{\bE\left[\left(b[x_{s}^{i, N},
\nu_{s}^{N}] - b[x_{s}^{\om_{i}}, P_{s}]\right)^{2}\right]\d s}\right.\\ &+
\left.\int_{0}^{t}{\bE\left[\left(c(x_{s}^{i, N}, \om_{i}) - c(x_{s}^{\om_{i}}, \om_{i})\right)^{2}\right]\d
s}\right),\\
	&\leq C \left(2\int_{0}^{t}{\bE\left[\sup_{u\leq s}\left|x_{u}^{i, N} -
x_{u}^{\om_{i}}\right|^{2}\right]\d s}\right.+ \left.\int_{0}^{t}{\sup_{1\leq j\leq N}\bE\left[\sup_{u\leq
s}\left|x_{u}^{\om_{j}} - x_{u}^{j, N}\right|^{2}\right]\d s}\right.\\&+
\left.\int_{0}^{t}{\bE\left[\left(b[x_{s}^{\om_{i}}, \nu_{s}] - b[x_{s}^{\om_{i}}, P_{s}]\right)^{2}\right]\d
s}\right).
\end{align*}
Applying Gronwall's Lemma to $\sup_{1\leq j\leq N}\bE\left[\sup_{u\leq
t}\left|x_{u}^{\om_{j}} - x_{u}^{j, N}\right|^{2}\right]$, it suffices to prove that for some $Z_{N}$:
\[\int_{0}^{t}{\bE\left[\left(b[x_{s}^{\om_{i}}, \nu_{s}] - b[x_{s}^{\om_{i}}, P_{s}]\right)^{2}\right]\d s}
\leq C/N + Z_{N}(\om_{1}, \dots, \om_{N}).\]
Indeed, for all $1\leq i\leq N$, (we write $x^{i}$ instead of $x^{\om_{i}}$ to simplify notations):
\[u_{i, N}:=\left(b[x_{s}^{i}, \nu_{s}] - b[x_{s}^{i}, P_{s}]\right)^{2} = \frac{1}{N^{2}}
\left(\Som{j}{1}{N}{T(x^{i}, x^{j})^{2}} + \sum_{k\neq l}{T(x^{i}, x^{k})T(x^{i}, x^{l})}\right),\]
where $T(x^{i}, x^{j}):=b(x^{i}_{s}, x^{j}_{s}, \om_{j}) - \int_{}{b(x^{i}_{s}, y, \pi) P_{s}(\d y, \d\pi)}$.
Since $b$ is bounded, we see that the first term is of order $(1/N)$. We only have to study the remaining
term:
	\[\bE\left[\underset{k\neq l}{\sum}T(x^{i}, x^{k})T(x^{i}, x^{l})\right]\leq CN +
\bE\left[\underset{\underset{k\neq l}{k\neq i, l\neq i}}{\sum}T(x^{i}, x^{k})T(x^{i}, x^{l})\right].\]
Since the $(x^{i})$ are independent, if we take conditional expectation w.r.t. $(x^{r}, r\neq l)$ in the last
term, we get:
\begin{align*}
	\bE\left[\underset{\underset{k\neq l}{k\neq i, l\neq i}}{\sum} {T(x^{i}, x^{k})T(x^{i},
x^{l})}\right]&= \bE\left[\bE\left[\left.\sum T(x^{i}, x^{k})T(x^{i}, x^{l})\right|x^{r}, r\neq
l\right]\right],\\
	&=\bE\left[\underset{\underset{k\neq l}{k\neq i,  l\neq i}}{\sum}T(x^{i}, x^{k})G_{l}(x^{i})\right] =
\bE\left[\underset{\underset{k\neq l}{k\neq i, l\neq i}}{\sum}G_{k}(x^{i})G_{l}(x^{i})\right],
\end{align*}
where $G_{l}(x) = G(x, \om_{l}) = \int_{}{b(x, y, \om_{l}) P_{s}^{\om_{l}}(\d y)} - \int_{}{b(x, y,
\pi)P_{s}(\d y,\d\pi)}$.
Defining \[Z_{N}(\om_{1}, \dots, \om_{N}):= \frac{C}{N}
\int_{0}^{T}{\bE\left[\left(\frac{1}{\sqrt{N}}\sum_{l=1}^{N} G(x^{i}_{s}, \om_{l})\right)^{2}\right]\d s},\]
in order to prove \eqref{eq:estimZN} it suffices to show that for some constant $C$,
\[\bbE\left[\int_{0}^{T}{\bE\left[\left(\frac{1}{\sqrt{N}}\Som{l}{1}{N}{G(x_{s}^{i},
\om_{l})}\right)^{2}\right]\d s}\right]\leq C.\]
The rest of the proof is devoted to prove this last assertion: we have successively (setting
$U_{N}(x^{\om_{i}}_{s}, \sou{\om}):= \frac{1}{\sqrt{N}}\Som{l}{1}{N}{G(x_{s}^{\om_{i}}, \om_{l})}$)
\begin{align*}
	\bbE\left[\int_{0}^{T}{ \bE\left[U_{N}(x_{s}^{\om_{i}}, \sou{\om})^{2}\right]\d s}\right]&\leq
\int_{0}^{T}{\bE\left[\bbE\left[U_{N}(x_{s}^{\om_{i}}, \sou{\om})^{2}\right]\right]\d s},\\
	&\leq \frac{1}{N} \int_{0}^{T}{\bE\left[\bbE\left[\Som{k}{1}{N}{\Som{l}{1}{N}{G(x_{s}^{\om_{i}},
\om_{k})G(x^{\om_{i}}_{s}, \om_{l})}}\right]\right]\d s},\\
	&\leq \	\frac{1}{N}\int_{0}^{T}{\bE\left[\bbE\left[\Som{l}{1}{N}{G(x^{\om_{i}}_{s},
\om_{l})^{2}}\right]\right]\d s} + C\\
	&+ \frac{1}{N} \int_{0}^{T}{\bE\left[\bbE\left[\sum_{l\neq k,\ l\neq i,\ k\neq i}{G(x_{s}^{\om_{i}},
\om_{k})G(x^{\om_{i}}_{s}, \om_{l})}\right]\right]\d s}.
\end{align*}
The first term of the RHS of the last inequality is bounded, since $b$ is bounded. But, if we condition w.r.t.
$\om_{r}$ for $r\neq i, r\neq k$, we see that the second term is zero. The result follows.

\section*{Acknowledgements} This is a part of my PhD thesis. I would like to thank my PhD supervisors
Giambattista Giacomin and Lorenzo Zambotti for introducing this subject, for their useful advice, and for
their encouragement. I would also like to thank the referee for useful suggestions and comments.

\bibliographystyle{abbrv}
\bibliography{../biblio1.bib}

\end{document}